\documentstyle[12pt]{amsart}
\numberwithin{equation}{section} 
\newtheorem{thm}[equation]{Theorem} 
\newtheorem{norithm}{Theorem (Nori)} 

\newtheorem{cor}[equation]{Corollary} 
\newtheorem{noricor}{Corollary (Nori)}

\newtheorem{lem}[equation]{Lemma}

\newtheorem{prop}[equation]{Proposition}

\theoremstyle{definition}

\theoremstyle{remark}

\newtheorem{rem}{Remark}    
\newtheorem{rems}{Remarks}

\begin{document}
\raggedbottom
\voffset=-.7truein
\hoffset=0truein
\vsize=8truein
\hsize=6truein
\textheight=8truein
\textwidth=6truein
\baselineskip=18truept


\def\C {\Bbb C} 
\def\Cn {\Bbb C^n} 
\def\R {\Bbb R}
\def\Rn {\Bbb R^n}
\def\Z {\Bbb Z}

\def\dbar {\bar \partial }
\def\dir {\cal D}
\def\lev#1{\cal L(#1)}
\def\lap {\Delta }
\def\ol {{\cal O}}
\def\E {{\cal E}}
\def\J {{\cal J}}
\def\U {{\cal U}}
\def\V {{\cal V}}
\def\z {\zeta } 
\def\Harm {\text {Harm}\, } 
\def\grad {\nabla }
\def\dexh {\{ M_k \} _{k=0}^{\infty } }
\def\sing#1{#1_{\text {sing}}}
\def\reg#1{#1_{\text {reg}}} 

\def\setof#1#2{\{ \, #1 \mid #2 \, \} }

\def\holecl {M\setminus \overline M_0}
\def\hole {M\setminus M_0}

\def\nd{\frac {\partial }{\partial\nu } }
\def\ndof#1{\frac {\partial#1}{\partial\nu } }

\def\pdof#1#2{\frac {\partial#1}{\partial#2}}

\def\cinf{C^{\infty }}

\def\diam{\text {diam} \, }

\def\real{\text {Re}\, }

\def\imag{\text {Im}\, }

\def\supp{\text {supp}\, }

\def\Vol{\text {\rm vol} \, }




\def\anal{analytic }
\def\analns{analytic}

\def\bdd{bounded }
\def\bddns{bounded}

\def\cpt{compact }
\def\cptns{compact}

\def\cpx{complex }
\def\cpxns{complex}

\def\cont{continuous }
\def\contns{continuous}

\def\dime{dimension }
\def\dimens{dimension }

\def\exh{exhaustion }
\def\exhns{exhaustion}

\def\fn{function }
\def\fnns{function}

\def\fns{functions }
\def\fnsns{functions}

\def\holo{holomorphic }
\def\holons{holomorphic}

\def\mero{meromorphic }
\def\merons{meromorphic}

\def\holoconvex{holomorphically convex }
\def\holoconvexns{holomorphically convex}

\def\ircomp{irreducible component }
\def\concomp{connected component }
\def\ircompns{irreducible component}
\def\concompns{connected component}
\def\ircomps{irreducible components }
\def\concomps{connected components }
\def\ircompsns{irreducible components}
\def\concompsns{connected components}

\def\irred{irreducible }
\def\irredns{irreducible}

\def\con{connected }
\def\conns{connected}

\def\comp{component }
\def\compns{component}
\def\comps{components }
\def\compsns{components}

\def\mfld{manifold }
\def\mfldns{manifold}
\def\mflds{manifolds }
\def\mfldsns{manifolds}

\def\nbd{neighborhood }
\def\nbds{neighborhoods }
\def\nbdns{neighborhood}
\def\nbdsns{neighborhoods}

\def\harm{harmonic }
\def\harmns{harmonic}
\def\plh{pluriharmonic }
\def\plhns{pluriharmonic}
\def\plsh{plurisubharmonic }
\def\plshns{plurisubharmonic}

\def\qplsh#1{$#1$-plurisubharmonic}
\def\hplsh{$(n-1)$-plurisubharmonic }
\def\hplshns{$(n-1)$-plurisubharmonic}

\def\para{parabolic }
\def\parans{parabolic}

\def\rel{relatively }
\def\relns{relatively}

\def\str{strictly }
\def\strns{strictly}

\def\wrt{with respect to }
\def\wrtns{with respect to}

\def\st {such that }
\def\stns {such that}

\def\hm {harmonic measure }
\def\hmns {harmonic measure}

\def\hmib {harmonic measure of the ideal boundary of }
\def\hmibns {harmonic measure of the ideal boundary of}


\def\atil{\tilde a}
\def\btil{\tilde b}
\def\ctil{\tilde c}
\def\dtil{\tilde d}
\def\etil{\tilde e}
\def\ftil{\tilde f}
\def\gtil{\tilde g}
\def\htil{\tilde h}
\def\itil{\tilde i}
\def\jtil{\tilde j}
\def\ktil{\tilde k}
\def\ltil{\tilde l}
\def\mtil{\tilde m}
\def\ntil{\tilde n}
\def\otil{\tilde o}
\def\ptil{\tilde p}
\def\qtil{\tilde q}
\def\rtil{\tilde r}
\def\stil{\tilde s}
\def\ttil{\tilde t}
\def\util{\tilde u}
\def\vtil{\tilde v}
\def\wtil{\tilde w}
\def\xtil{\tilde x}
\def\ytil{\tilde y}
\def\ztil{\tilde z}

\def\Atil{\tilde A}
\def\Btil{\widetilde B}
\def\Ctil{\widetilde C}
\def\Dtil{\widetilde D}
\def\Etil{\widetilde E}
\def\Ftil{\widetilde F}
\def\Gtil{\widetilde G}
\def\Htil{\widetilde H}
\def\Itil{\tilde I}
\def\Jtil{\widetilde J}
\def\Ktil{\widetilde K}
\def\Ltil{\widetilde L}
\def\Mtil{\widetilde M}
\def\Ntil{\widetilde N}
\def\Otil{\widetilde O}
\def\Ptil{\widetilde P}
\def\Qtil{\widetilde Q}
\def\Rtil{\widetilde R}
\def\Stil{\widetilde S}
\def\Ttil{\widetilde T}
\def\Util{\widetilde U}
\def\Vtil{\widetilde V}
\def\Wtil{\widetilde W}
\def\Xtil{\widetilde X}
\def\Ytil{\widetilde Y}
\def\Ztil{\widetilde Z}

\def\alphatil {\tilde \alpha  }
\def\betatil {\tilde \beta  }
\def\gammatil {\tilde \gamma  }
\def\deltatil {\tilde \delta }
\def\epsilontil {\tilde \epsilon  }
\def\varepsilontil {\tilde \varepsilon  }
\def\zetatil {\tilde \zeta  }
\def\etatil {\tilde \eta  }
\def\thetatil {\tilde \theta  }
\def\varthetatil {\tilde \vartheta  }
\def\iotatil {\tilde \iota  }
\def\kappatil {\tilde \kappa  }
\def\lambdatil {\tilde \lambda  }
\def\mutil {\tilde \mu  }
\def\nutil {\tilde \nu  }
\def\xitil {\tilde \xi  }
\def\pitil {\tilde \pi  }
\def\varpitil {\tilde \varpi  }
\def\rhotil {\tilde \rho  }
\def\varrhotil {\tilde \varrho  }
\def\sigmatil {\tilde \sigma  }
\def\varsigmatil {\tilde \varsigma  }
\def\tautil {\tilde \tau  }
\def\upsilontil {\tilde \upsilon  }
\def\phitil {\tilde \phi  }
\def\varphitil {\tilde \varphi  }
\def\chitil {\tilde \chi  }
\def\psitil {\tilde \psi  }
\def\omegatil {\tilde \omega  }

\def\Gammatil {\widetilde \Gamma  }
\def\Deltatil {\tilde \Delta }
\def\Thetatil {\widetilde \Theta  }
\def\Lambdatil {\tilde \Lambda  }
\def\Xitil {\widetilde \Xi  }
\def\Pitil {\widetilde \Pi  }
\def\Sigmatil {\widetilde \Sigma  }
\def\Upsilontil {\widetilde \Upsilon  }
\def\Phitil {\tilde \Phi  }
\def\Psitil {\widetilde \Psi  }
\def\Omegatil {\widetilde \Omega  }

\def\varGammatil {\widetilde \varGamma  }
\def\varDeltatil {\tilde \varDelta  }
\def\varThetatil {\widetilde \varTheta  }
\def\varLambdatil {\tilde \varLambda  }
\def\varXitil {\widetilde \varXi  }
\def\varPitil {\widetilde \varPi  }
\def\varSigmatil {\widetilde \varSigma  }
\def\varUpsilontil {\widetilde \varUpsilon  }
\def\varPhitil {\tilde \varPhi  }
\def\varPsitil {\widetilde \varPsi  }
\def\varOmegatil {\widetilde \varOmega  }

\def\boldGammatil {\widetilde \boldGamma  }
\def\boldDeltatil {\tilde \boldDelta  }
\def\boldThetatil {\widetilde \boldTheta  }
\def\boldLambdatil {\tilde \boldLambda  }
\def\boldXitil {\widetilde \boldXi  }
\def\boldPitil {\widetilde \boldPi  }
\def\boldSigmatil {\widetilde \boldSigma  }
\def\boldUpsilontil {\widetilde \boldUpsilon  }
\def\boldPhitil {\tilde \boldPhi  }
\def\boldPsitil {\widetilde \boldPsi  }
\def\boldOmegatil {\widetilde \boldOmega  }


\def\ahat{\hat a}
\def\bhat{\hat b}
\def\chat{\hat c}
\def\dhat{\hat d}
\def\ehat{\hat e}
\def\fhat{\hat f}
\def\ghat{\hat g}
\def\hhat{\hat h}
\def\ihat{\hat i}
\def\jhat{\hat j}
\def\khat{\hat k}
\def\lhat{\hat l}
\def\mhat{\hat m}
\def\nhat{\hat n}
\def\ohat{\hat o}
\def\phat{\hat p}
\def\qhat{\hat q}
\def\rhat{\hat r}
\def\shat{\hat s}
\def\that{\hat t}
\def\uhat{\hat u}
\def\vhat{\hat v}
\def\what{\hat w}
\def\xhat{\hat x}
\def\yhat{\hat y}
\def\zhat{\hat z}

\def\Ahat{\hat A}
\def\Bhat{\widehat B}
\def\Chat{\widehat C}
\def\Dhat{\widehat D}
\def\Ehat{\widehat E}
\def\Fhat{\widehat F}
\def\Ghat{\widehat G}
\def\Hhat{\widehat H}
\def\Ihat{\hat I}
\def\Jhat{\widehat J}
\def\Khat{\widehat K}
\def\Lhat{\widehat L}
\def\Mhat{\widehat M}
\def\Nhat{\widehat N}
\def\Ohat{\widehat O}
\def\Phat{\widehat P}
\def\Qhat{\widehat Q}
\def\Rhat{\widehat R}
\def\Shat{\widehat S}
\def\That{\widehat T}
\def\Uhat{\widehat U}
\def\Vhat{\widehat V}
\def\What{\widehat W}
\def\Xhat{\widehat X}
\def\Yhat{\widehat Y}
\def\Zhat{\widehat Z}

\def\alphahat {\hat \alpha  }
\def\betahat {\hat \beta  }
\def\gammahat {\hat \gamma  }
\def\deltahat {\hat \delta }
\def\epsilonhat {\hat \epsilon  }
\def\varepsilonhat {\hat \varepsilon  }
\def\zetahat {\hat \zeta  }
\def\etahat {\hat \eta  }
\def\thetahat {\hat \theta  }
\def\varthetahat {\hat \vartheta  }
\def\iotahat {\hat \iota  }
\def\kappahat {\hat \kappa  }
\def\lambdahat {\hat \lambda  }
\def\muhat {\hat \mu  }
\def\nuhat {\hat \nu  }
\def\xihat {\hat \xi  }
\def\pihat {\hat \pi  }
\def\varpihat {\hat \varpi  }
\def\rhohat {\hat \rho  }
\def\varrhohat {\hat \varrho  }
\def\sigmahat {\hat \sigma  }
\def\varsigmahat {\hat \varsigma  }
\def\tauhat {\hat \tau  }
\def\upsilonhat {\hat \upsilon  }
\def\phihat {\hat \phi  }
\def\varphihat {\hat \varphi  }
\def\vphihat {\hat \varphi  }
\def\chihat {\hat \chi  }
\def\psihat {\hat \psi  }
\def\omegahat {\hat \omega  }

\def\Gammahat {\widehat \Gamma  }
\def\Deltahat {\hat \Delta }
\def\Thetahat {\widehat \Theta  }
\def\Lambdahat {\hat \Lambda  }
\def\Xihat {\widehat \Xi  }
\def\Pihat {\widehat \Pi  }
\def\Sigmahat {\widehat \Sigma  }
\def\Upsilonhat {\widehat \Upsilon  }
\def\Phihat {\hat \Phi  }
\def\Psihat {\widehat \Psi  }
\def\Omegahat {\widehat \Omega  }

\def\varGammahat {\widehat \varGamma  }
\def\varDeltahat {\hat \varDelta  }
\def\varThetahat {\widehat \varTheta  }
\def\varLambdahat {\hat \varLambda  }
\def\varXihat {\widehat \varXi  }
\def\varPihat {\widehat \varPi  }
\def\varSigmahat {\widehat \varSigma  }
\def\varUpsilonhat {\widehat \varUpsilon  }
\def\varPhihat {\hat \varPhi  }
\def\varPsihat {\widehat \varPsi  }
\def\varOmegahat {\widehat \varOmega  }

\def\boldGammahat {\widehat \boldGamma  }
\def\boldDeltahat {\hat \boldDelta  }
\def\boldThetahat {\widehat \boldTheta  }
\def\boldLambdahat {\hat \boldLambda  }
\def\boldXihat {\widehat \boldXi  }
\def\boldPihat {\widehat \boldPi  }
\def\boldSigmahat {\widehat \boldSigma  }
\def\boldUpsilonhat {\widehat \boldUpsilon  }
\def\boldPhihat {\hat \boldPhi  }
\def\boldPsihat {\widehat \boldPsi  }
\def\boldOmegahat {\widehat \boldOmega  }

\def\seq #1#2{\{#1_{#2}\} }

\def\aseq{ \{ a_{\nu } \} }
\def\bseq{ \{ b_{\nu } \} }
\def\cseq{ \{ c_{\nu } \} }
\def\dseq{ \{ d_{\nu } \} }
\def\eseq{ \{ e_{\nu } \} }
\def\fseq{ \{ f_{\nu } \} }
\def\gseq{ \{ g_{\nu } \} }
\def\hseq{ \{ h_{\nu } \} }
\def\iseq{ \{ i_{\nu } \} }
\def\jseq{ \{ j_{\nu } \} }
\def\kseq{ \{ k_{\nu } \} }
\def\lseq{ \{ l_{\nu } \} }
\def\mseq{ \{ m_{\nu } \} }
\def\nseq{ \{ n_{\nu } \} }
\def\oseq{ \{ o_{\nu } \} }
\def\pseq{ \{ p_{\nu } \} }
\def\qseq{ \{ q_{\nu } \} }
\def\rseq{ \{ r_{\nu } \} }
\def\sseq{ \{ s_{\nu } \} }
\def\tseq{ \{ t_{\nu } \} }
\def\useq{ \{ u_{\nu } \} }
\def\vseq{ \{ v_{\nu } \} }
\def\wseq{ \{ w_{\nu } \} }
\def\xseq{ \{ x_{\nu } \} }
\def\yseq{ \{ y_{\nu } \} }
\def\zseq{ \{ z_{\nu } \} }

\def\Aseq{ \{ A_{\nu } \} }
\def\Bseq{ \{ B_{\nu } \} }
\def\Cseq{ \{ C_{\nu } \} }
\def\Dseq{ \{ D_{\nu } \} }
\def\Eseq{ \{ E_{\nu } \} }
\def\Fseq{ \{ F_{\nu } \} }
\def\Gseq{ \{ G_{\nu } \} }
\def\Hseq{ \{ H_{\nu } \} }
\def\Iseq{ \{ I_{\nu } \} }
\def\Jseq{ \{ J_{\nu } \} }
\def\Kseq{ \{ K_{\nu } \} }
\def\Lseq{ \{ L_{\nu } \} }
\def\Mseq{ \{ M_{\nu } \} }
\def\Nseq{ \{ N_{\nu } \} }
\def\Oseq{ \{ O_{\nu } \} }
\def\Pseq{ \{ P_{\nu } \} }
\def\Qseq{ \{ Q_{\nu } \} }
\def\Rseq{ \{ R_{\nu } \} }
\def\Sseq{ \{ S_{\nu } \} }
\def\Tseq{ \{ T_{\nu } \} }
\def\Useq{ \{ U_{\nu } \} }
\def\Vseq{ \{ V_{\nu } \} }
\def\Wseq{ \{ W_{\nu } \} }
\def\Xseq{ \{ X_{\nu } \} }
\def\Yseq{ \{ Y_{\nu } \} }
\def\Zseq{ \{ Z_{\nu } \} }

\def\alphaseq { \{ \alpha _{\nu } \} }
\def\betaseq { \{ \beta _{\nu } \} }
\def\gammaseq { \{ \gamma _{\nu } \} }
\def\deltaseq { \{ \delta _{\nu } \} }
\def\epsilonseq { \{ \epsilon _{\nu } \} }
\def\varepsilonseq { \{ \varepsilon _{\nu } \} }
\def\zetaseq { \{ \zeta _{\nu } \} }
\def\etaseq { \{ \eta _{\nu } \} }
\def\thetaseq { \{ \theta _{\nu } \} }
\def\varthetaseq { \{ \vartheta _{\nu } \} }
\def\iotaseq { \{ \iota _{\nu } \} }
\def\kappaseq { \{ \kappa _{\nu } \} }
\def\lambdaseq { \{ \lambda _{\nu } \} }
\def\museq { \{ \mu _{\nu } \} }
\def\nuseq { \{ \nu  _{\nu } \} }
\def\xiseq { \{ \xi _{\nu } \} }
\def\piseq { \{ \pi _{\nu } \} }
\def\varpiseq { \{ \varpi _{\nu } \} }
\def\rhoseq { \{ \rho _{\nu } \} }
\def\varrhoseq { \{ \varrho _{\nu } \} }
\def\sigmaseq { \{ \sigma _{\nu } \} }
\def\varsigmaseq { \{ \varsigma _{\nu } \} }
\def\tauseq { \{ \tau _{\nu } \} }
\def\upsilonseq { \{ \upsilon _{\nu } \} }
\def\phiseq { \{ \phi _{\nu } \} }
\def\varphiseq { \{ \varphi _{\nu } \} }
\def\chiseq { \{ \chi _{\nu } \} }
\def\psiseq { \{ \psi _{\nu } \} }
\def\omegaseq { \{ \omega _{\nu } \} }

\def\Gammaseq { \{ \Gamma _{\nu } \} }
\def\Deltaseq { \{ \Delta _{\nu } \} }
\def\Thetaseq { \{ \Theta _{\nu } \} }
\def\Lambdaseq { \{ \Lambda _{\nu } \} }
\def\Xiseq { \{ \Xi _{\nu } \} }
\def\Piseq { \{ \Pi _{\nu } \} }
\def\Sigmaseq { \{ \Sigma _{\nu } \} }
\def\Upsilonseq { \{ \Upsilon _{\nu } \} }
\def\Phiseq { \{ \Phi _{\nu } \} }
\def\Psiseq { \{ \Psi _{\nu } \} }
\def\Omegaseq { \{ \Omega _{\nu } \} }

\def\varGammaseq { \{ \varGamma _{\nu } \} }
\def\varDeltaseq { \{ \varDelta _{\nu } \} }
\def\varThetaseq { \{ \varTheta _{\nu } \} }
\def\varLambdaseq { \{ \varLambda _{\nu } \} }
\def\varXiseq { \{ \varXi _{\nu } \} }
\def\varPiseq { \{ \varPi _{\nu } \} }
\def\varSigmaseq { \{ \varSigma _{\nu } \} }
\def\varUpsilonseq { \{ \varUpsilon _{\nu } \} }
\def\varPhiseq { \{ \varPhi _{\nu } \} }
\def\varPsiseq { \{ \varPsi _{\nu } \} }
\def\varOmegaseq { \{ \varOmega _{\nu } \} }

\def\boldGammaseq { \{ \boldGamma _{\nu } \} }
\def\boldDeltaseq { \{ \boldDelta _{\nu } \} }
\def\boldThetaseq { \{ \boldTheta _{\nu } \} }
\def\boldLambdaseq { \{ \boldLambda _{\nu } \} }
\def\boldXiseq { \{ \boldXi _{\nu } \} }
\def\boldPiseq { \{ \boldPi _{\nu } \} }
\def\boldSigmaseq { \{ \boldSigma _{\nu } \} }
\def\boldUpsilonseq { \{ \boldUpsilon _{\nu } \} }
\def\boldPhiseq { \{ \boldPhi _{\nu } \} }
\def\boldPsiseq { \{ \boldPsi _{\nu } \} }
\def\boldOmegaseq { \{ \boldOmega _{\nu } \} }

\def\amu{   a_{\mu }  }
\def\bmu{   b_{\mu }  }
\def\cmu{   c_{\mu }  }
\def\dmu{   d_{\mu }  }
\def\emu{   e_{\mu }  }
\def\fmu{   f_{\mu }  }
\def\gmu{   g_{\mu }  }
\def\hmu{   h_{\mu }  }
\def\imu{   i_{\mu }  }
\def\jmu{   j_{\mu }  }
\def\kmu{   k_{\mu }  }
\def\lmu{   l_{\mu }  }
\def\mmu{   m_{\mu }  }
\def\nmu{   n_{\mu }  }
\def\omu{   o_{\mu }  }
\def\pmu{   p_{\mu }  }
\def\qmu{   q_{\mu }  }
\def\rmu{   r_{\mu }  }
\def\smu{   s_{\mu }  }
\def\tmu{   t_{\mu }  }
\def\umu{   u_{\mu }  }
\def\vmu{   v_{\mu }  }
\def\wmu{   w_{\mu }  }
\def\xmu{   x_{\mu }  }
\def\ymu{   y_{\mu }  }
\def\zmu{   z_{\mu }  }

\def\Amu{   A_{\mu }  }
\def\Bmu{   B_{\mu }  }
\def\Cmu{   C_{\mu }  }
\def\Dmu{   D_{\mu }  }
\def\Emu{   E_{\mu }  }
\def\Fmu{   F_{\mu }  }
\def\Gmu{   G_{\mu }  }
\def\Hmu{   H_{\mu }  }
\def\Imu{   I_{\mu }  }
\def\Jmu{   J_{\mu }  }
\def\Kmu{   K_{\mu }  }
\def\Lmu{   L_{\mu }  }
\def\Mmu{   M_{\mu }  }
\def\Nmu{   N_{\mu }  }
\def\Omu{   O_{\mu }  }
\def\Pmu{   P_{\mu }  }
\def\Qmu{   Q_{\mu }  }
\def\Rmu{   R_{\mu }  }
\def\Smu{   S_{\mu }  }
\def\Tmu{   T_{\mu }  }
\def\Umu{   U_{\mu }  }
\def\Vmu{   V_{\mu }  }
\def\Wmu{   W_{\mu }  }
\def\Xmu{   X_{\mu }  }
\def\Ymu{   Y_{\mu }  }
\def\Zmu{   Z_{\mu }  }

\def\alphamu{\alpha _{\mu }}
\def\betamu{\beta _{\mu }}
\def\gammamu{\gamma _{\mu }}
\def\deltamu{\delta _{\mu }}
\def\epsilonmu{\epsilon _{\mu }}
\def\varepsilonmu{\varepsilon _{\mu }}
\def\zetamu{\zeta _{\mu }}
\def\etamu{\eta _{\mu }}
\def\thetamu{\theta _{\mu }}
\def\varthetamu{\vartheta _{\mu }}
\def\iotamu{\iota _{\mu }}
\def\kappamu{\kappa _{\mu }}
\def\lambdamu{\lambda _{\mu }}
\def\mumu{\mu _{\mu }}
\def\numu{\nu _{\mu }}
\def\ximu{\xi _{\mu }}
\def\pimu{\pi _{\mu }}
\def\varpimu{\varpi _{\mu }}
\def\rhomu{\rho _{\mu }}
\def\varrhomu{\varrho _{\mu }}
\def\sigmamu{\sigma _{\mu }}
\def\varsigmamu{\varsigma _{\mu }}
\def\taumu{\tau _{\mu }}
\def\upsilonmu{\upsilon _{\mu }}
\def\phimu{\phi _{\mu }}
\def\varphimu{\varphi _{\mu }}
\def\chimu{\chi _{\mu }}
\def\psimu{\psi _{\mu }}
\def\omegamu{\omega _{\mu }}

\def\Gammamu{\Gamma _{\mu }}
\def\Deltamu{\Delta _{\mu }}
\def\Thetamu{\Theta _{\mu }}
\def\Lambdamu{\Lambda _{\mu }}
\def\Ximu{\Xi _{\mu }}
\def\Pimu{\Pi _{\mu }}
\def\Sigmamu{\Sigma _{\mu }}
\def\Upsilonmu{\Upsilon _{\mu }}
\def\Phimu{\Phi _{\mu }}
\def\Psimu{\Psi _{\mu }}
\def\Omegamu{\Omega _{\mu }}

\def\varGammamu{\varGamma _{\mu }}
\def\varDeltamu{\varDelta _{\mu }}
\def\varThetamu{\varTheta _{\mu }}
\def\varLambdamu{\varLambda _{\mu }}
\def\varXimu{\varXi _{\mu }}
\def\varPimu{\varPi _{\mu }}
\def\varSigmamu{\varSigma _{\mu }}
\def\varUpsilonmu{\varUpsilon _{\mu }}
\def\varPhimu{\varPhi _{\mu }}
\def\varPsimu{\varPsi _{\mu }}
\def\varOmegamu{\varOmega _{\mu }}

\def\boldGammamu{\boldGamma _{\mu }}
\def\boldDeltamu{\boldDelta _{\mu }}
\def\boldThetamu{\boldTheta _{\mu }}
\def\boldLambdamu{\boldLambda _{\mu }}
\def\boldXimu{\boldXi _{\mu }}
\def\boldPimu{\boldPi _{\mu }}
\def\boldSigmamu{\boldSigma _{\mu }}
\def\boldUpsilonmu{\boldUpsilon _{\mu }}
\def\boldPhimu{\boldPhi _{\mu }}
\def\boldPsimu{\boldPsi _{\mu }}
\def\boldOmegamu{\boldOmega _{\mu }}


\def\asmu{   a^{(\mu )}  }
\def\bsmu{   b^{(\mu )}  }
\def\csmu{   c^{(\mu )}  }
\def\dsmu{   d^{(\mu )}  }
\def\esmu{   e^{(\mu )}  }
\def\fsmu{   f^{(\mu )}  }
\def\gsmu{   g^{(\mu )}  }
\def\hsmu{   h^{(\mu )}  }
\def\ismu{   i^{(\mu )}  }
\def\jsmu{   j^{(\mu )}  }
\def\ksmu{   k^{(\mu )}  }
\def\lsmu{   l^{(\mu )}  }
\def\msmu{   m^{(\mu )}  }
\def\nsmu{   n^{(\mu )}  }
\def\osmu{   o^{(\mu )}  }
\def\psmu{   p^{(\mu )}  }
\def\qsmu{   q^{(\mu )}  }
\def\rsmu{   r^{(\mu )}  }
\def\ssmu{   s^{(\mu )}  }
\def\tsmu{   t^{(\mu )}  }
\def\usmu{   u^{(\mu )}  }
\def\vsmu{   v^{(\mu )}  }
\def\wsmu{   w^{(\mu )}  }
\def\xsmu{   x^{(\mu )}  }
\def\ysmu{   y^{(\mu )}  }
\def\zsmu{   z^{(\mu )}  }

\def\Asmu{   A^{(\mu )}  }
\def\Bsmu{   B^{(\mu )}  }
\def\Csmu{   C^{(\mu )}  }
\def\Dsmu{   D^{(\mu )}  }
\def\Esmu{   E^{(\mu )}  }
\def\Fsmu{   F^{(\mu )}  }
\def\Gsmu{   G^{(\mu )}  }
\def\Hsmu{   H^{(\mu )}  }
\def\Ismu{   I^{(\mu )}  }
\def\Jsmu{   J^{(\mu )}  }
\def\Ksmu{   K^{(\mu )}  }
\def\Lsmu{   L^{(\mu )}  }
\def\Msmu{   M^{(\mu )}  }
\def\Nsmu{   N^{(\mu )}  }
\def\Osmu{   O^{(\mu )}  }
\def\Psmu{   P^{(\mu )}  }
\def\Qsmu{   Q^{(\mu )}  }
\def\Rsmu{   R^{(\mu )}  }
\def\Ssmu{   S^{(\mu )}  }
\def\Tsmu{   T^{(\mu )}  }
\def\Usmu{   U^{(\mu )}  }
\def\Vsmu{   V^{(\mu )}  }
\def\Wsmu{   W^{(\mu )}  }
\def\Xsmu{   X^{(\mu )}  }
\def\Ysmu{   Y^{(\mu )}  }
\def\Zsmu{   Z^{(\mu )}  }

\def\alphasmu{\alpha ^{(\mu )}}
\def\betasmu{\beta ^{(\mu )}}
\def\gammasmu{\gamma ^{(\mu )}}
\def\deltasmu{\delta ^{(\mu )}}
\def\epsilonsmu{\epsilon ^{(\mu )}}
\def\varepsilonsmu{\varepsilon ^{(\mu )}}
\def\zetasmu{\zeta ^{(\mu )}}
\def\etasmu{\eta ^{(\mu )}}
\def\thetasmu{\theta ^{(\mu )}}
\def\varthetasmu{\vartheta ^{(\mu )}}
\def\iotasmu{\iota ^{(\mu )}}
\def\kappasmu{\kappa ^{(\mu )}}
\def\lambdasmu{\lambda ^{(\mu )}}
\def\musmu{\mu ^{(\mu )}}
\def\nusmu{\nu ^{(\mu )}}
\def\xismu{\xi ^{(\mu )}}
\def\pismu{\pi ^{(\mu )}}
\def\varpismu{\varpi ^{(\mu )}}
\def\rhosmu{\rho ^{(\mu )}}
\def\varrhosmu{\varrho ^{(\mu )}}
\def\sigmasmu{\sigma ^{(\mu )}}
\def\varsigmasmu{\varsigma ^{(\mu )}}
\def\tausmu{\tau ^{(\mu )}}
\def\upsilonsmu{\upsilon ^{(\mu )}}
\def\phismu{\phi ^{(\mu )}}
\def\varphismu{\varphi ^{(\mu )}}
\def\chismu{\chi ^{(\mu )}}
\def\psismu{\psi ^{(\mu )}}
\def\omegasmu{\omega ^{(\mu )}}

\def\Gammasmu{\Gamma ^{(\mu )}}
\def\Deltasmu{\Delta ^{(\mu )}}
\def\Thetasmu{\Theta ^{(\mu )}}
\def\Lambdasmu{\Lambda ^{(\mu )}}
\def\Xismu{\Xi ^{(\mu )}}
\def\Pismu{\Pi ^{(\mu )}}
\def\Sigmasmu{\Sigma ^{(\mu )}}
\def\Upsilonsmu{\Upsilon ^{(\mu )}}
\def\Phismu{\Phi ^{(\mu )}}
\def\Psismu{\Psi ^{(\mu )}}
\def\Omegasmu{\Omega ^{(\mu )}}

\def\varGammasmu{\varGamma ^{(\mu )}}
\def\varDeltasmu{\varDelta ^{(\mu )}}
\def\varThetasmu{\varTheta ^{(\mu )}}
\def\varLambdasmu{\varLambda ^{(\mu )}}
\def\varXismu{\varXi ^{(\mu )}}
\def\varPismu{\varPi ^{(\mu )}}
\def\varSigmasmu{\varSigma ^{(\mu )}}
\def\varUpsilonsmu{\varUpsilon ^{(\mu )}}
\def\varPhismu{\varPhi ^{(\mu )}}
\def\varPsismu{\varPsi ^{(\mu )}}
\def\varOmegasmu{\varOmega ^{(\mu )}}

\def\boldGammasmu{\boldGamma ^{(\mu )}}
\def\boldDeltasmu{\boldDelta ^{(\mu )}}
\def\boldThetasmu{\boldTheta ^{(\mu )}}
\def\boldLambdasmu{\boldLambda ^{(\mu )}}
\def\boldXismu{\boldXi ^{(\mu )}}
\def\boldPismu{\boldPi ^{(\mu )}}
\def\boldSigmasmu{\boldSigma ^{(\mu )}}
\def\boldUpsilonsmu{\boldUpsilon ^{(\mu )}}
\def\boldPhismu{\boldPhi ^{(\mu )}}
\def\boldPsismu{\boldPsi ^{(\mu )}}
\def\boldOmegasmu{\boldOmega ^{(\mu )}}

\def\anu{   a_{\nu }  }
\def\bnu{   b_{\nu }  }
\def\cnu{   c_{\nu }  }
\def\dnu{   d_{\nu }  }
\def\enu{   e_{\nu }  }
\def\fnu{   f_{\nu }  }
\def\gnu{   g_{\nu }  }
\def\hnu{   h_{\nu }  }
\def\inu{   i_{\nu }  }
\def\jnu{   j_{\nu }  }
\def\knu{   k_{\nu }  }
\def\lnu{   l_{\nu }  }
\def\mnu{   m_{\nu }  }
\def\nnu{   n_{\nu }  }
\def\onu{   o_{\nu }  }
\def\pnu{   p_{\nu }  }
\def\qnu{   q_{\nu }  }
\def\rnu{   r_{\nu }  }
\def\snu{   s_{\nu }  }
\def\tnu{   t_{\nu }  }
\def\unu{   u_{\nu }  }
\def\vnu{   v_{\nu }  }
\def\wnu{   w_{\nu }  }
\def\xnu{   x_{\nu }  }
\def\ynu{   y_{\nu }  }
\def\znu{   z_{\nu }  }

\def\Anu{   A_{\nu }  }
\def\Bnu{   B_{\nu }  }
\def\Cnu{   C_{\nu }  }
\def\Dnu{   D_{\nu }  }
\def\Enu{   E_{\nu }  }
\def\Fnu{   F_{\nu }  }
\def\Gnu{   G_{\nu }  }
\def\Hnu{   H_{\nu }  }
\def\Inu{   I_{\nu }  }
\def\Jnu{   J_{\nu }  }
\def\Knu{   K_{\nu }  }
\def\Lnu{   L_{\nu }  }
\def\Mnu{   M_{\nu }  }
\def\Nnu{   N_{\nu }  }
\def\Onu{   O_{\nu }  }
\def\Pnu{   P_{\nu }  }
\def\Qnu{   Q_{\nu }  }
\def\Rnu{   R_{\nu }  }
\def\Snu{   S_{\nu }  }
\def\Tnu{   T_{\nu }  }
\def\Unu{   U_{\nu }  }
\def\Vnu{   V_{\nu }  }
\def\Wnu{   W_{\nu }  }
\def\Xnu{   X_{\nu }  }
\def\Ynu{   Y_{\nu }  }
\def\Znu{   Z_{\nu }  }

\def\alphanu{\alpha _{\nu }}
\def\betanu{\beta _{\nu }}
\def\gammanu{\gamma _{\nu }}
\def\deltanu{\delta _{\nu }}
\def\epsilonnu{\epsilon _{\nu }}
\def\varepsilonnu{\varepsilon _{\nu }}
\def\zetanu{\zeta _{\nu }}
\def\etanu{\eta _{\nu }}
\def\thetanu{\theta _{\nu }}
\def\varthetanu{\vartheta _{\nu }}
\def\iotanu{\iota _{\nu }}
\def\kappanu{\kappa _{\nu }}
\def\lambdanu{\lambda _{\nu }}
\def\munu{\mu _{\nu }}
\def\nunu{\nu _{\nu }}
\def\xinu{\xi _{\nu }}
\def\pinu{\pi _{\nu }}
\def\varpinu{\varpi _{\nu }}
\def\rhonu{\rho _{\nu }}
\def\varrhonu{\varrho _{\nu }}
\def\sigmanu{\sigma _{\nu }}
\def\varsigmanu{\varsigma _{\nu }}
\def\taunu{\tau _{\nu }}
\def\upsilonnu{\upsilon _{\nu }}
\def\phinu{\phi _{\nu }}
\def\varphinu{\varphi _{\nu }}
\def\chinu{\chi _{\nu }}
\def\psinu{\psi _{\nu }}
\def\omeganu{\omega _{\nu }}

\def\Gammanu{\Gamma _{\nu }}
\def\Deltanu{\Delta _{\nu }}
\def\Thetanu{\Theta _{\nu }}
\def\Lambdanu{\Lambda _{\nu }}
\def\Xinu{\Xi _{\nu }}
\def\Pinu{\Pi _{\nu }}
\def\Sigmanu{\Sigma _{\nu }}
\def\Upsilonnu{\Upsilon _{\nu }}
\def\Phinu{\Phi _{\nu }}
\def\Psinu{\Psi _{\nu }}
\def\Omeganu{\Omega _{\nu }}

\def\varGammanu{\varGamma _{\nu }}
\def\varDeltanu{\varDelta _{\nu }}
\def\varThetanu{\varTheta _{\nu }}
\def\varLambdanu{\varLambda _{\nu }}
\def\varXinu{\varXi _{\nu }}
\def\varPinu{\varPi _{\nu }}
\def\varSigmanu{\varSigma _{\nu }}
\def\varUpsilonnu{\varUpsilon _{\nu }}
\def\varPhinu{\varPhi _{\nu }}
\def\varPsinu{\varPsi _{\nu }}
\def\varOmeganu{\varOmega _{\nu }}

\def\boldGammanu{\boldGamma _{\nu }}
\def\boldDeltanu{\boldDelta _{\nu }}
\def\boldThetanu{\boldTheta _{\nu }}
\def\boldLambdanu{\boldLambda _{\nu }}
\def\boldXinu{\boldXi _{\nu }}
\def\boldPinu{\boldPi _{\nu }}
\def\boldSigmanu{\boldSigma _{\nu }}
\def\boldUpsilonnu{\boldUpsilon _{\nu }}
\def\boldPhinu{\boldPhi _{\nu }}
\def\boldPsinu{\boldPsi _{\nu }}
\def\boldOmeganu{\boldOmega _{\nu }}


\def\asnu{   a^{(\nu )}  }
\def\bsnu{   b^{(\nu )}  }
\def\csnu{   c^{(\nu )}  }
\def\dsnu{   d^{(\nu )}  }
\def\esnu{   e^{(\nu )}  }
\def\fsnu{   f^{(\nu )}  }
\def\gsnu{   g^{(\nu )}  }
\def\hsnu{   h^{(\nu )}  }
\def\isnu{   i^{(\nu )}  }
\def\jsnu{   j^{(\nu )}  }
\def\ksnu{   k^{(\nu )}  }
\def\lsnu{   l^{(\nu )}  }
\def\msnu{   m^{(\nu )}  }
\def\nsnu{   n^{(\nu )}  }
\def\osnu{   o^{(\nu )}  }
\def\psnu{   p^{(\nu )}  }
\def\qsnu{   q^{(\nu )}  }
\def\rsnu{   r^{(\nu )}  }
\def\ssnu{   s^{(\nu )}  }
\def\tsnu{   t^{(\nu )}  }
\def\usnu{   u^{(\nu )}  }
\def\vsnu{   v^{(\nu )}  }
\def\wsnu{   w^{(\nu )}  }
\def\xsnu{   x^{(\nu )}  }
\def\ysnu{   y^{(\nu )}  }
\def\zsnu{   z^{(\nu )}  }

\def\Asnu{   A^{(\nu )}  }
\def\Bsnu{   B^{(\nu )}  }
\def\Csnu{   C^{(\nu )}  }
\def\Dsnu{   D^{(\nu )}  }
\def\Esnu{   E^{(\nu )}  }
\def\Fsnu{   F^{(\nu )}  }
\def\Gsnu{   G^{(\nu )}  }
\def\Hsnu{   H^{(\nu )}  }
\def\Isnu{   I^{(\nu )}  }
\def\Jsnu{   J^{(\nu )}  }
\def\Ksnu{   K^{(\nu )}  }
\def\Lsnu{   L^{(\nu )}  }
\def\Msnu{   M^{(\nu )}  }
\def\Nsnu{   N^{(\nu )}  }
\def\Osnu{   O^{(\nu )}  }
\def\Psnu{   P^{(\nu )}  }
\def\Qsnu{   Q^{(\nu )}  }
\def\Rsnu{   R^{(\nu )}  }
\def\Ssnu{   S^{(\nu )}  }
\def\Tsnu{   T^{(\nu )}  }
\def\Usnu{   U^{(\nu )}  }
\def\Vsnu{   V^{(\nu )}  }
\def\Wsnu{   W^{(\nu )}  }
\def\Xsnu{   X^{(\nu )}  }
\def\Ysnu{   Y^{(\nu )}  }
\def\Zsnu{   Z^{(\nu )}  }

\def\alphasnu{\alpha ^{(\nu )}}
\def\betasnu{\beta ^{(\nu )}}
\def\gammasnu{\gamma ^{(\nu )}}
\def\deltasnu{\delta ^{(\nu )}}
\def\epsilonsnu{\epsilon ^{(\nu )}}
\def\varepsilonsnu{\varepsilon ^{(\nu )}}
\def\zetasnu{\zeta ^{(\nu )}}
\def\etasnu{\eta ^{(\nu )}}
\def\thetasnu{\theta ^{(\nu )}}
\def\varthetasnu{\vartheta ^{(\nu )}}
\def\iotasnu{\iota ^{(\nu )}}
\def\kappasnu{\kappa ^{(\nu )}}
\def\lambdasnu{\lambda ^{(\nu )}}
\def\musnu{\mu ^{(\nu )}}
\def\nusnu{\nu ^{(\nu )}}
\def\xisnu{\xi ^{(\nu )}}
\def\pisnu{\pi ^{(\nu )}}
\def\varpisnu{\varpi ^{(\nu )}}
\def\rhosnu{\rho ^{(\nu )}}
\def\varrhosnu{\varrho ^{(\nu )}}
\def\sigmasnu{\sigma ^{(\nu )}}
\def\varsigmasnu{\varsigma ^{(\nu )}}
\def\tausnu{\tau ^{(\nu )}}
\def\upsilonsnu{\upsilon ^{(\nu )}}
\def\phisnu{\phi ^{(\nu )}}
\def\varphisnu{\varphi ^{(\nu )}}
\def\chisnu{\chi ^{(\nu )}}
\def\psisnu{\psi ^{(\nu )}}
\def\omegasnu{\omega ^{(\nu )}}

\def\Gammasnu{\Gamma ^{(\nu )}}
\def\Deltasnu{\Delta ^{(\nu )}}
\def\Thetasnu{\Theta ^{(\nu )}}
\def\Lambdasnu{\Lambda ^{(\nu )}}
\def\Xisnu{\Xi ^{(\nu )}}
\def\Pisnu{\Pi ^{(\nu )}}
\def\Sigmasnu{\Sigma ^{(\nu )}}
\def\Upsilonsnu{\Upsilon ^{(\nu )}}
\def\Phisnu{\Phi ^{(\nu )}}
\def\Psisnu{\Psi ^{(\nu )}}
\def\Omegasnu{\Omega ^{(\nu )}}

\def\varGammasnu{\varGamma ^{(\nu )}}
\def\varDeltasnu{\varDelta ^{(\nu )}}
\def\varThetasnu{\varTheta ^{(\nu )}}
\def\varLambdasnu{\varLambda ^{(\nu )}}
\def\varXisnu{\varXi ^{(\nu )}}
\def\varPisnu{\varPi ^{(\nu )}}
\def\varSigmasnu{\varSigma ^{(\nu )}}
\def\varUpsilonsnu{\varUpsilon ^{(\nu )}}
\def\varPhisnu{\varPhi ^{(\nu )}}
\def\varPsisnu{\varPsi ^{(\nu )}}
\def\varOmegasnu{\varOmega ^{(\nu )}}

\def\boldGammasnu{\boldGamma ^{(\nu )}}
\def\boldDeltasnu{\boldDelta ^{(\nu )}}
\def\boldThetasnu{\boldTheta ^{(\nu )}}
\def\boldLambdasnu{\boldLambda ^{(\nu )}}
\def\boldXisnu{\boldXi ^{(\nu )}}
\def\boldPisnu{\boldPi ^{(\nu )}}
\def\boldSigmasnu{\boldSigma ^{(\nu )}}
\def\boldUpsilonsnu{\boldUpsilon ^{(\nu )}}
\def\boldPhisnu{\boldPhi ^{(\nu )}}
\def\boldPsisnu{\boldPsi ^{(\nu )}}
\def\boldOmegasnu{\boldOmega ^{(\nu )}}


\def\vphi {\varphi }


\def\inv{   ^{-1}  }

\def\ainv{   a^{-1}  }
\def\binv{   b^{-1}  }
\def\cinv{   c^{-1}  }
\def\dinv{   d^{-1}  }
\def\einv{   e^{-1}  }
\def\finv{   f^{-1}  }
\def\ginv{   g^{-1}  }
\def\hinv{   h^{-1}  }
\def\iinv{   i^{-1}  }
\def\jinv{   j^{-1}  }
\def\kinv{   k^{-1}  }
\def\linv{   l^{-1}  }
\def\minv{   m^{-1}  }
\def\ninv{   n^{-1}  }
\def\oinv{   o^{-1}  }
\def\pinv{   p^{-1}  }
\def\qinv{   q^{-1}  }
\def\rinv{   r^{-1}  }
\def\sinv{   s^{-1}  }
\def\tinv{   t^{-1}  }
\def\uinv{   u^{-1}  }
\def\vinv{   v^{-1}  }
\def\winv{   w^{-1}  }
\def\xinv{   x^{-1}  }
\def\yinv{   y^{-1}  }
\def\zinv{   z^{-1}  }

\def\Ainv{   A^{-1}  }
\def\Binv{   B^{-1}  }
\def\Cinv{   C^{-1}  }
\def\Dinv{   D^{-1}  }
\def\Einv{   E^{-1}  }


\def\Ginv{   G^{-1}  }
\def\Hinv{   H^{-1}  }
\def\Iinv{   I^{-1}  }
\def\Jinv{   J^{-1}  }
\def\Kinv{   K^{-1}  }
\def\Linv{   L^{-1}  }
\def\Minv{   M^{-1}  }
\def\Ninv{   N^{-1}  }
\def\Oinv{   O^{-1}  }
\def\Pinv{   P^{-1}  }
\def\Qinv{   Q^{-1}  }
\def\Rinv{   R^{-1}  }
\def\Sinv{   S^{-1}  }
\def\Tinv{   T^{-1}  }
\def\Uinv{   U^{-1}  }
\def\Vinv{   V^{-1}  }
\def\Winv{   W^{-1}  }
\def\Xinv{   X^{-1}  }
\def\Yinv{   Y^{-1}  }
\def\Zinv{   Z^{-1}  }

\def\alphainv{\alpha ^{-1}}
\def\betainv{\beta ^{-1}}
\def\gammainv{\gamma ^{-1}}
\def\deltainv{\delta ^{-1}}
\def\epsiloninv{\epsilon ^{-1}}
\def\varepsiloninv{\varepsilon ^{-1}}
\def\zetainv{\zeta ^{-1}}
\def\etainv{\eta ^{-1}}
\def\thetainv{\theta ^{-1}}
\def\varthetainv{\vartheta ^{-1}}
\def\iotainv{\iota ^{-1}}
\def\kappainv{\kappa ^{-1}}
\def\lambdainv{\lambda ^{-1}}
\def\muinv{\mu ^{-1}}
\def\nuinv{\nu ^{-1}}
\def\xiinv{\xi ^{-1}}
\def\piinv{\pi ^{-1}}
\def\varpiinv{\varpi ^{-1}}
\def\rhoinv{\rho ^{-1}}
\def\varrhoinv{\varrho ^{-1}}
\def\sigmainv{\sigma ^{-1}}
\def\varsigmainv{\varsigma ^{-1}}
\def\tauinv{\tau ^{-1}}
\def\upsiloninv{\upsilon ^{-1}}
\def\phiinv{\phi ^{-1}}
\def\varphiinv{\varphi ^{-1}}
\def\vphiinv{\varphi ^{-1}}
\def\chiinv{\chi ^{-1}}
\def\psiinv{\psi ^{-1}}
\def\omegainv{\omega ^{-1}}

\def\Gammainv{\Gamma ^{-1}}
\def\Deltainv{\Delta ^{-1}}
\def\Thetainv{\Theta ^{-1}}
\def\Lambdainv{\Lambda ^{-1}}
\def\Xiinv{\Xi ^{-1}}
\def\Piinv{\Pi ^{-1}}
\def\Sigmainv{\Sigma ^{-1}}
\def\Upsiloninv{\Upsilon ^{-1}}
\def\Phiinv{\Phi ^{-1}}
\def\Psiinv{\Psi ^{-1}}
\def\Omegainv{\Omega ^{-1}}

\def\varGammainv{\varGamma ^{-1}}
\def\varDeltainv{\varDelta ^{-1}}
\def\varThetainv{\varTheta ^{-1}}
\def\varLambdainv{\varLambda ^{-1}}
\def\varXiinv{\varXi ^{-1}}
\def\varPiinv{\varPi ^{-1}}
\def\varSigmainv{\varSigma ^{-1}}
\def\varUpsiloninv{\varUpsilon ^{-1}}
\def\varPhiinv{\varPhi ^{-1}}
\def\varPsiinv{\varPsi ^{-1}}
\def\varOmegainv{\varOmega ^{-1}}

\def\boldGammainv{\boldGamma ^{-1}}
\def\boldDeltainv{\boldDelta ^{-1}}
\def\boldThetainv{\boldTheta ^{-1}}
\def\boldLambdainv{\boldLambda ^{-1}}
\def\boldXiinv{\boldXi ^{-1}}
\def\boldPiinv{\boldPi ^{-1}}
\def\boldSigmainv{\boldSigma ^{-1}}
\def\boldUpsiloninv{\boldUpsilon ^{-1}}
\def\boldPhiinv{\boldPhi ^{-1}}
\def\boldPsiinv{\boldPsi ^{-1}}
\def\boldOmegainv{\boldOmega ^{-1}}
\title[weak Lefschetz theorems]
{The $L^2$\ $\dbar $-method, weak Lefschetz theorems,
and the topology of K\"ahler manifolds}
\author[T.~Napier]{Terrence Napier}
\address{Department of Mathematics\\Lehigh University\\Bethlehem, PA 18015}
\email{tjn2@@lehigh.edu}
\thanks{Research partially 
supported by NSF grant DMS9411154}
\author[M.~Ramachandran]{Mohan Ramachandran}
\address{Department of Mathematics\\SUNY at Buffalo\\Buffalo, NY 14214}
\email{ramac-m@@newton.math.buffalo.edu}
\thanks{Research partially 
supported by NSF grant DMS9626169}


\date{}

\maketitle

\begin{abstract}
A new approach to Nori's
weak Lefschetz theorem is described. The new approach,
which involves the $\dbar $-method, avoids moving
arguments and gives much stronger results. 
In particular, it is proved that if $X$ and $Y$ are \con smooth 
projective varieties of positive dimension
and 
$f : Y @>>> X$ is a \holo immersion with ample normal bundle, 
then the image of $\pi_1(Y)$ in $\pi _1(X)$ is of finite index.
This result is obtained as a consequence of a direct 
generalization of Nori's theorem. 
The second part concerns a new approach to
the theorem of Burns which states that a 
quotient of the unit ball in 
$\C ^n$ ($n\geq 3$)
by a discrete group of automorphisms which has a strongly
pseudoconvex boundary component has only finitely many ends.
The following generalization is obtained. 
If a complete Hermitian manifold~$X$ of dimension
$n\geq 3$ has a strongly pseudoconvex end~$E$ and 
$\text {Ricci}\, (X) \leq -C$ for some positive constant~$C$, 
then, away from~$E$, $X$ has finite 
volume. 
\end{abstract}

\section{Introduction} \label{intro}

In~[No], Nori studied the fundamental group of complements of nodal
curves with ample normal bundle in smooth projective surfaces. 
The main tool was the following weak Lefschetz theorem:
\begin{norithm}
Suppose $\Phi : U @>>> X$
is a local biholomorphism from a \con \cpx manifold~$U$ into a \con 
smooth projective variety~$X$ of dimension at least~$2$ and $U$
contains a \con effective divisor~$Y$ with \cpt support and ample normal 
bundle. Then, for every Zariski open subset~$Z$
of~$X$, the image of $\pi _1(\Phi\inv (Z))$ in $\pi _1(Z)$ is of finite
index. 
\end{norithm}
For $X$ a surface, he obtained sharp bounds for the index
using the Hodge index theorem. 
A striking corollary of this result is the following:
\begin{noricor}
If $X$ and $Y$ are \con smooth projective
varieties with 
$$
\dim X=\dim Y +1>1
$$ 
and 
$f : Y @>>> X$ is a \holo immersion with ample normal bundle, 
then the image of $\pi_1(Y)$ in $\pi _1(X)$ is of finite index. 
\end{noricor}

Nori's proof of these results depends heavily on deformations. 
The first step 
is to show that a large multiple of the divisor~$Y$ in the theorem moves
in a family in which the general member is irreducible and meets~$Y$
and the union of these members contains an open subset of~$U$.
Unfortunately, moving arguments do not seem to apply in the 
higher codimensional case, because Fulton and Lazarsfeld~[FL2]
have observed that for a certain smooth projective $4$-fold and a smooth 
surface~$Y$
in~$X$ with ample normal bundle constructed by Gieseker,
no multiple of~$Y$ in~$X$ moves. 
Given the existence of a sufficiently large number of
deformations, the rest of the proof of Nori's
weak Lefschetz theorem has been streamlined by 
Campana~[C1] and Koll\'ar~[K]. In~[NR], 
another proof of Nori's theorem
was given when $Z=X$ using harmonic \fnsns ,
but it was the same in spirit as the earlier arguments. 
A survey on Lefschetz type theorems can be found in Fulton~[F].

In this paper we introduce a new approach which avoids moving arguments and
which gives much stronger results. In particular, the new approach allows 
one to address the case of higher codimension. Before giving precise 
statements, we recall some terminology. Let $Y$ be a \cpx \anal 
subspace of complex space~$U$. We denote the structure sheaf of~$U$
by $\ol _U$ and the ideal sheaf of~$Y$ in~$U$ by~$\cal I_Y$. The 
{\it formal completion}~$\Uhat $~{\it of}~$U$ {\it \wrtns }~$Y$ 
is the ringed space
$$
(\Uhat , \ol _{\Uhat })=(Y, \lim _{ @<<< }\ol _U/\cal I_Y^n).
$$
If $\cal F$ is an \anal sheaf on~$U$ we denote by $\widehat {\cal F}$ 
the associated \anal sheaf on~$\Uhat $ given by
$$
\widehat {\cal F} = 
\lim _{ @<<< }(\cal F \otimes \ol _U/\cal I_Y^n).
$$
If $\cal F$ is coherent, then $\widehat {\cal F}$ is also coherent
over~$\ol _{\Uhat }$. The main
result is the following generalization of Nori's weak Lefschetz theorem:
\begin{thm}
\it 
Suppose $\Phi : U @>>> X$
is a \holo map from a \con \cpx manifold~$U$ into a \con 
smooth projective variety~$X$ of dimension at least~$2$ 
which is a submersion at some point. Let~$Y\subset U$ be a
\con \cpt \anal subspace \st 
$\dim H^0(\Uhat , \widehat {\cal L}) <\infty $
for every locally free \anal sheaf~$\cal L$ on~$U$. 
Then, for every Zariski open subset~$Z$
of~$X$, the image of $\pi _1(\Phi\inv (Z))$ in $\pi _1(Z)$ is of finite
index. 
\end{thm}

\begin{rems} 1. For example, by a theorem of Hartshorne~[H]
(and Grothendieck~[Gr]),
$H^0(\Uhat , \widehat {\cal L})$ is finite dimensional 
when $Y$ is a \con \cpt \anal subspace which is locally a complete
intersection and which has ample normal bundle
(or even $k$-ample normal bundle
in the sense of Sommese~[So] where $k=\dim Y-1$).

\noindent 2. Theorem~0.1 also holds for $U$ \irred and reduced
and $X$ normal and projective. Moreover, as will be shown in 
Sect.~3 (Corollary~3.4),  in the smooth case one only 
needs finite dimensionality for $\cal L$ the \anal pullback
of an  invertible sheaf on~$X$. 

\noindent 3. As a consequence of Theorem~0.1, one can remove the dimension
restriction on the subspace~$Y$ in the corollary to Nori's theorem.
More precisely, we get the following:
\end{rems} 
\begin{cor}
If $X$ and $Y$ are \con smooth projective
varieties of positive dimension and 
$f : Y @>>> X$ is a \holo immersion with ample normal bundle, 
then the image of $\pi_1(Y)$ in $\pi _1(X)$ is of finite index. 
\end{cor}
Hironaka and Matsumura~[HM] proved
the analogous result
for algebraic fundamental groups when $f$ is an inclusion
with ample normal bundle. However, the result 
for topological fundamental groups (as stated in the
above corollary)
is new (provided $\dim X > \dim Y +1$)
even for $f$ an inclusion. 
Moreover, simple examples show that,
if $\dim X > \dim Y +1$, then, even if $f$ is an inclusion
(with ample normal bundle), the 
map $\pi_1(Y) @>>> \pi _1(X)$ is not necessarily surjective.

The idea of the proof of Theorem~0.1 is to form
a covering space~$\Ztil @>>> Z$ with fundamental group equal to
the image~$G$ of $\pi _1(\Phi\inv (Z))$ and then
construct $L^2$ \holo sections of a 
suitable line bundle which separate the sheets of the covering. 
This construction is a standard application of the 
$L^2$~$\dbar $-method 
(Andreotti-Vesentini~[AV], 
H\"ormander~[Ho], Skoda~[Sk], Demailly~[D1]).
Pulling these sections back to $\Phi\inv (Z)$ by a lifting of~$\Phi $, 
the finite dimensionality of the space of \holo
sections on the formal completion
gives a bound on the dimension of the space of sections on~$\Ztil $
and hence a bound on the degree of the covering space
(i.e.~on the index of~$G$). 

\begin{rem}
Campana~[C2] has independently applied $L^2$-methods
to study exceptional curves on coverings of surfaces. 
\end{rem}

The second main result of this paper generalizes a theorem of 
Burns~[B]
which states that a quotient of the unit ball in~$\C ^n$ ($n\geq 3$)
by a discrete group of automorphisms which has a strongly pseudoconvex
boundary component has only finitely many ends. The main tools are
a theorem of Lempert on the compactification of a pseudoconvex 
boundary from the pseudoconcave side~[L], a finiteness theorem
of Andreotti for pseudoconcave manifolds~[A], and the $L^2$ Riemann-Roch 
inequality of Nadel and Tsuji~[NT]. The precise statement is 
as follows:
\begin{thm}
If a complete Hermitian manifold~$X$ of (complex) dimension
at least~$3$ has a strongly pseudoconvex end and 
$\text {Ricci}\, (X) \leq -C$ for some positive constant~$C$, then, away
from the strongly pseudoconvex end, the manifold has finite volume.
\end{thm}
 

As in the proof of Theorem~0.1, the idea is to apply finite dimensionality
of the space of \holo sections of a line bundle. By Lempert's theorem, 
one can cap off the strongly pseudoconvex end by a domain in a smooth
projective variety.  Andreotti's finiteness theorem applied to the 
resulting pseudoconcave manifold gives finite dimensionality of the space
of \holo sections of a suitable line bundle. Finally, the $L^2$ Riemann-Roch 
inequality of Nadel and Tsuji gives a (finite) upper bound for the volume
in terms of the dimension of this space of sections.

\begin{rem}
One natural question which arises is might 
there be
an improved version of the $L^2$ Riemann-Roch inequality which would
give improved bounds for the volume in Theorem~0.3
as well as the index in~Theorem~0.1? Also, 
for $X$ a surface in the corollary to Nori's weak Lefschetz 
theorem, 
Nori~[No] found bounds for the index in terms of certain intersection
numbers. It is therefore 
natural to look for analogous bounds in more general cases.
\end{rem}

Sect.~1 begins with a proof of Theorem~0.1 in the case where~$\Phi $
is a local biholomorphism. The main idea of the 
new approach is easy to
see in this context, and, although a few technicalities 
arise in the general case, the proof is essentially the same. 
The proof of Theorem~0.1 is then given. 
The required result from the $L^2$ \  $\dbar $-method is discussed 
in Sect.~2. Further generalizations of the weak Lefschetz theorem
for $X$ not necessarily projective are considered in Sect.~3.
Theorem~0.3 is proved in Sect.~4, which may be read independently
of Sects.~1--3.

\noindent {\it Acknowledgements}. 
Madhav Nori suggested we reformulate
Theorem~0.1 in terms of formal completions, which considerably widened 
its scope.
Charles Epstein
told us about Lempert's result. For this and other useful advice, we would
like to thank them both.  We would also like to thank Alan Nadel for 
bringing the $L^2$ Riemann-Roch inequality to our attention, 
Dan Burns for useful discussions on his theorem, and
Raghavan Narasimhan for his interest in this work. 
Finally, we would like to thank the referee for helpful
suggestions.

\section{Weak Lefschetz theorems for a projective variety}

This section contains the proof Theorem~0.1. We first 
prove the theorem for the case of a local biholomorphism. This is a direct 
generalization 
to immersed \cpx spaces of {\it arbitrary} codimension
(Nori proved the theorem stated below for $Y$ an ample divisor in~$U$).
More general versions 
will be stated later. Aside
from a few minor technical problems, however, the proofs of all of
the generalizations are the same in spirit as the 
proof of this special case.

\begin{thm}
Let $U$ be a \con \cpx manifold, let~$X$ be 
a \con smooth projective variety of dimension~$n>1$, let  
$\Phi : U @>>> X$ be a \holo map, let~$Y$ be a \con \cpt \anal subspace
(not necessarily reduced) of~$U$, and let $\Uhat $ be the formal completion
of~$U$ with respect to~$Y$.  Assume that 
\begin{enumerate}
\item[(i)] $\Phi $ is locally bi\holons , and
\item[(ii)] $\dim H^0(\Uhat , \widehat {\ol (\Phi ^*L)})<\infty $
for every \holo line bundle~$L$ on~$X$.  
\end{enumerate} 
Then there is a positive constant~$b$ depending only on the mapping
$\Phi :U @>>> X$ and the subspace~$Y\subset U$ \stns , if $R\subset X$
is a nowhere dense analytic subset of~$X$ and~$V$ is a \con \nbd of~$Y$ in~$U$,
then the image~$G$ of 
$\pi _1(V\setminus \Phi\inv (R)) @>>> \pi _1(X\setminus R)$
is of index at most~$b$ in~$\pi _1(X\setminus R)$. Moreover, if
$\Phi (Y)\cap R=\emptyset $, then the image of 
$\pi _1(Y) @>>> \pi _1(X\setminus R)$ is also of index at most~$b$. 
\end{thm}

\begin{pf}
Given $R$,$V$, and~$G$ as in the statement of the theorem,
let $S=\Phi\inv (R)$, 
let $M=X\setminus R$, let $W=V\setminus S$,
and let $\pi : \Mtil @>>> M$ be a \con covering space
with $\pi _* (\pi _1(\Mtil ))=G$.  Thus $\pi : \Mtil @>>> M$ has
degree $d=[\pi _1(M):G]$ and we have the
following commutative diagram of \holo mappings: 
\begin{center}\begin{picture}(250,80)
\put(5,10){$W=V\setminus S\subset V $}
\put(90,14){\vector(1,0){50}}
\put(145,10){$X \supset X\setminus R=M$}
\put(110,3){$\Phi $}
\put(125,55){\vector(3,-1){90}}
\put(110,60){$\Mtil $}
\put(17,25){\vector(3,1){90}}
\put(57,45){$\Phitil $}
\put(170,45){$\pi $}
\end{picture}
\end{center}
%
%
%
Since $X$ is projective,
there exists a Hermitian \holo line bundle~$(L,h)$ with
positive curvature and a K\"ahler metric~$g$ on~$X$. 
As will be shown in Sect.~2 (see Corollary~2.3),
the $L^2$~$\dbar $-method, in the form given by
Skoda~[Sk] and Demailly~[D1], enables one to prove that there is 
a positive integer~$\nu $ independent 
of~$R$~and~$V$ such that
$$
d\leq \dim  H^0_{L^2}(\Mtil , \ol (\pi ^* (L^{\nu }\otimes K_M)));
$$
where $K_M$ is the canonical bundle on~$M$ and 
$H^0_{L^2}(\Mtil , \ol (\pi ^* (L^{\nu }\otimes K_M)))$ is the space
of \holo sections of $\pi ^* (L^{\nu }\otimes K_M)$ which are 
in~$L^2$ with respect to the Hermitian metrics $\pi ^*(h\otimes g^*)$
on $\pi ^* (L^{\nu }\otimes K_M)$ and~$\pi ^* g$ on~$\Mtil $. 
If
$s\in H^0_{L^2}(\Mtil , \ol (\pi ^* (L^{\nu }\otimes K_M)))$, then
$\Phitil ^* s$ is a \holo section of $\Phi ^* (L^{\nu }\otimes K_X)$
on $W$. Given a point $x_0\in S\cap V$, $\Phi $ maps a \nbdns~$Q$
of~$x_0$ in~$V$ biholomorphically onto $\Phi (Q) \subset X$. Hence $\Phitil $
maps $Q\setminus S$ biholomorphically onto its image in~$\Mtil $ and,
therefore, $\Phitil ^*s$ is in~$L^2$ on $Q\setminus S$ with respect
to the Hermitian metrics $\Phi ^*(h\otimes g^*)$ in 
$\Phi ^* (L^{\nu }\otimes K_X)$ and $\Phi ^*g$ on~$U$. Since these metrics
are defined over the entire set~$U$ and a square integrable function
which is \holo outside a nowhere dense  analytic set in a manifold extends
holomorphically past the analytic set, $\Phitil ^*s$ 
extends to a \holo section of $\Phi ^* (L^{\nu }\otimes K_X)$ on~$V$. 
Therefore
$$
d\leq \dim H^0(V, \ol (\Phi ^* (L^{\nu }\otimes K_X))).
$$
On the other hand, 
by a general fact
about formal completions, if $\cal F$ is a coherent \anal sheaf on~$V$, 
then the kernel of the mapping
$$
H^0(V,\cal F) @>>> H^0(\Vhat ,\widehat {\cal F})=H^0(\Uhat ,\widehat {\cal F})
$$
consists of all of the sections of~$\cal F$ on~$V$ which vanish 
on a \nbd of~$Y$ in~$V$ (see [BS, 
Proposition VI.2.7]). In particular, if~$\cal F$ is locally 
free, then this mapping is injective. 
Therefore, taking
$\cal F=\ol (\Phi ^*(L^{\nu }\otimes K_X))$, we get
$$
d\leq \dim H^0(V,\cal F) \leq \dim  
H^0(\Uhat ,\widehat {\cal F})<\infty .
$$
Thus $b=\dim  H^0(\Uhat ,\widehat {\cal F})$ is a uniform bound
for~$d$ independent of~$R$ and~$V$.

Finally, if $\Phi (Y)\cap R=\emptyset $, then we may choose the \nbdns~$V$
so that $V\subset U\setminus S$ and the map $\pi _1(Y) @>>> \pi _1(V)$
is a surjective isomorphism. Hence the image of $\pi _1(Y)$ in
$\pi _1(M)$ is equal to the image of $\pi _1(V)=\pi _1(V\setminus S)$ and
therefore is of index at most~$b$. 
\end{pf}

We now consider generalizations.
If in the above theorem one assumes only that~$\Phi $ is a generic
submersion (or a generic local biholomorphism), then a slight technical problem
arises.  While (as one may easily check) the section $\Phitil ^*s$ 
of $\Phi ^*(L^\nu \otimes K_X)$ extends holomorphically past~$S$
near points at which~$\Phi $ is submersive, $\Phitil ^*s$  
need not extend near points where $\text {rank} \, \Phi _* < n$. 
However, as we will see, $\Phitil ^*s$ does  
extend as a \holo $n$-form with values in~$\Phi ^*L^\nu $.
A simple illustration is given by 
$$
U=\Delta \ni z \overset {\Phi } \mapsto \zeta =z^2 \in \Delta , \quad 
\Mtil =\Delta ^* \ni z \overset {\pi } \mapsto \zeta =z^2 \in \Delta ^*=M,
\quad \text {and } s=z\inv \pi ^* d\zeta .
$$
Here, $\Phitil ^*s$  does not extend as a section of the pullback
of the canonical bundle, but the corresponding \holo $1$-form~$2dz$ 
does extend. 

In fact, by passing to desingularizations, one also gets this
extension property for $U$ and $X$ singular. 
Given an \irred reduced 
\cpx space~$A$ and a positive integer~$n$, we denote by~$\Omega ^n_A$
the coherent \anal sheaf on~$A$ obtained by forming a desingularization
$\check A @>>> A$ of~$A$ and taking the direct image of $\Omega ^n_{\check A}$.
By the following lemma, this sheaf is independent of the choice of the 
desingularization.
\begin{lem}[Grauert and Riemenschneider 
[GR, Sect.~2.1{]}]
Let $A$ be an \irred reduced \cpx space of dimension~$m$
and let~$n$ be a positive integer.
Suppose that, for $i=1,2$, $B_i$ is a 
\con \cpx manifold of dimension~$m$ and 
$\Psi _i : B_i @>>> A$ is a proper modification. 
Then 
$(\Psi _1)_*\Omega ^n_{B_1}=(\Psi _2)_*\Omega ^n_{B_2}$. 
\end{lem}
The proof is similar to the proof for $\dim A=n$ given 
in~[GR]. The main point is that if $A$
is smooth, then $\Psi _1$ is bi\holo outside an analytic
set of codimension at least~$2$ in~$A$. For the general
case, one passes to a common proper modification of 
$B_1$~and~$B_2$.  

We may now state the extension property as follows:
\begin{lem} 
\it 
Let $\Phi :U @>>> X$ be a \holo mapping of \irred 
reduced \cpx spaces $U$~and~$X$ of dimensions $m$~and~$n$, respectively,
\st $\Phi (U)$ has nonempty interior. Suppose
\begin{center}\begin{picture}(250,80)
\put(5,10){$W=U\setminus S\subset U$}
\put(90,14){\vector(1,0){50}}
\put(145,10){$X \supset X\setminus R=M$}
\put(110,3){$\Phi $}
\put(125,55){\vector(3,-1){90}}
\put(110,60){$\Mtil $}
\put(17,25){\vector(3,1){90}}
\put(57,45){$\Phitil $}
\put(170,45){$\pi $}
\end{picture}
\end{center}
%
%
is a commutative diagram of \holo mappings where
$R\subset X$ is a nowhere dense  \anal subset which contains~$\sing X$,  
$S=\Phi\inv (R)$, and 
$\pi : \Mtil @>>> M$ is a \con \holo covering space. 
Let $L$ be a \holo line bundle on~$X$ and let $\theta $ be a \holo 
$n$-form with values in $\pi^*L$ on $\Mtil $ which is in $L^2$
with respect to the liftings 
of a Hermitian metric~$h$ in~$L$ on~$X$ and a Hermitian metric~$g$ on~$M$. 
Then the pullback 
$(\Phitil | _{\reg W})^*\theta$ of~$\theta$ to a \holo
$n$-form with values in $\Phi ^*L$ on $\reg W$ extends to a (unique)
section in $H^0(U, \ol (\Phi ^*L) \otimes \Omega ^n_U)$. 
\end{lem}
The proof uses standard methods but will be postponed until
the end of this section (see also Sakai~[S]). We may now apply the 
argument given in the proof 
of Theorem~1.1 to get Theorem~0.1 of the introduction.
In fact, we get the following:
\begin{thm} Let $U$
be an \irred reduced \cpx space, let $X$ be a \con normal projective variety
of dimension~$n>1$,   
let $\Phi :U @>>> X$ be a \holo map, let~$Y$ be a \con \cpt \anal 
subspace (not necessarily reduced) of~$U$, and
let $\Uhat $ be the formal completion of~$U$ with respect to~$Y$.  
Assume that  
\begin{enumerate}
\item[(i)] $\Phi (U)$ has nonempty interior, and 
\item[(ii)] 
$\dim H^0(\Uhat , \widehat {\ol (\Phi ^*L)}\otimes \widehat {\Omega ^n_U})
<\infty $ for every \holo line bundle~$L$ on~$X$.
\end{enumerate}
Then there exists a positive
constant~$b$ depending only on the mapping $\Phi : U @>>> X$ and the 
subspace~$Y$ \stns , 
if~$R\subset X$ is a nowhere dense \anal subset and 
$V$ is a \con \nbd of~$Y$ in~$U$, 
then the image of 
$\pi _1(V\setminus \Phi\inv (R)) @>>> \pi _1 (X\setminus R)$
is of index at most~$b$. Moreover, if
$\Phi (Y)\cap R=\emptyset $, then the image of~$\pi _1(Y)$ in 
$\pi _1(X\setminus R)$ is also of index at most~$b$.
\end{thm}
%
%
\begin{pf} Given $R$ and $V$ as in the statement of the
theorem, we get a commutative diagram  
\begin{center}\begin{picture}(250,80)
\put(5,10){$W=V\setminus S\subset V$}
\put(90,14){\vector(1,0){50}}
\put(145,10){$X \supset X\setminus R=M$}
\put(110,3){$\Phi $}
\put(125,55){\vector(3,-1){90}}
\put(110,60){$\Mtil $}
\put(17,25){\vector(3,1){90}}
\put(57,45){$\Phitil $}
\put(170,45){$\pi $}
\end{picture}
\end{center}
as in the proof of Theorem~1.1.  
Since $X$ is normal, the map
$\pi _1(M\setminus \sing X) @>>> \pi _1(M)$ is surjective. 
Therefore, by 
replacing~$R$ by $R\cup\sing X$, we may assume that $\sing X\subset R$;
i.e.~that $M$ is a complete K\"ahler manifold. 
If 
$s\in H^0_{L^2}(\Mtil , \ol ((\pi ^*L^\nu )\otimes K_{\Mtil }))$ 
for some~$\nu $ (with respect to metrics lifted 
from the base), then, by Lemma~1.3, the pullback to $\reg W$ 
{\it as a \holo $n$-form
with values in~$\Phi ^*L^\nu $}
extends to a 
unique section in $H^0(V, \ol (\Phi ^*L^\nu ) \otimes \Omega ^n_U)$. 
By applying Corollary~2.3 as in the proof of Theorem~1.1, one now
gets the required bound on the index. 
\end{pf}

A finiteness theorem of  
Hartshorne~[H, Theorem~III.4.1]
and Grothendieck~[Gr] and the above theorem together 
imply immediately that, 
in Nori's weak Lefschetz theorem, one may take the mapping to
be a generic submersion and the
subvariety to be of arbitrary codimension. More precisely, we have
the following:
\begin{cor}
Let $U$ be a \con \cpx manifold, let~$X$ be 
a \con normal projective variety of dimension~$n>1$, let  
$\Phi : U @>>> X$ be a \holo map, and let~$Y$ be a 
positive dimensional \con \cpt \anal subspace
(not necessarily reduced) of~$U$. Assume that 
\begin{enumerate}
\item[(i)] $\Phi (U)$ has nonempty interior, 
\item[(ii)] $Y$ is locally a complete intersection in~$U$, and
\item[(iii)] The normal bundle $N_{Y/U}$ is ample.   
\end{enumerate} 
Then there is a positive constant~$b$ depending only on the mapping
$\Phi :U @>>> X$ and the subspace~$Y\subset U$ \stns , 
if~$Z$ is a nonempty Zariski open subset of~$X$ and 
$V$ is a \con \nbd of~$Y$ in~$U$, 
then the image of 
$\pi _1(V\cap \Phi\inv (Z)) @>>> \pi _1 (Z)$
is of index at most~$b$ in~$\pi _1(Z)$. Moreover, if
$\Phi (Y)\subset Z$, then the image of~$\pi _1(Y)$ in 
$\pi _1(Z)$ is also of index at most~$b$.
\end{cor}

\begin{rems}
1. The approach of considering sections of 
vector bundles
on formal completions
fits well with Grothendieck's approach to the Lefschetz theorems~[Gr]
(see also~[H]). In a sense, the results of this paper extend 
to the topological fundamental group 
Grothendieck's Lefschetz theorems concerning the algebraic fundamental group. 

\noindent 2. Further generalizations in which $X$ is not necessarily 
projective will be stated and proved in Sect.~3. A slightly more precise 
bound for the index in terms of the dimension of a 
space of sections will also be obtained. 
\end{rems}

We conclude this section with the proof of the extension property.
\begin{pf*}{Proof of Lemma~1.3}
We first observe that we may assume that 
$U$~and~$X$ are smooth and that $R$ is a divisor with normal
crossings by passing to desingularizations. More precisely, 
we may form a commutative diagram 
\begin{center}\begin{picture}(250,90)
\put(60,10){$U$}
\put(90,14){\vector(1,0){50}}
\put(150,10){$X$}
\put(115,0){$\Phi $}
\put(65,55){\vector(0,-1){30}}
\put(50,60){$U\times _XX'$}
\put(45,35){$\text{pr}_U$}
\put(100,64){\vector(1,0){40}}
\put(150,60){$X'$}
\put(115,75){$\text {pr}_{X'}$}
\put(160,35){$\beta $}
\put(155,55){\vector(0,-1){30}}
\end{picture}
\end{center}
%
%
where $X'$ is a \con \cpx manifold, $R'=\beta\inv (R)$ of $R$
is a divisor with normal crossings, and 
$\beta  : X' @>>> X$ is a proper modification which maps 
$M'=X'\setminus R'$ biholomorphically onto $M=X\setminus R$. 
Since
$\text {pr}_U \inv (W)=W\times _MM'$ is just the graph of the restriction
of~$\Phi $ to a mapping $W @>>> M'=M$ and $U$ is \irred,
$\text {pr}_U \inv (W)$ is an open \irred subset 
of $U\times _XX'$ which is mapped isomorphically onto~$W$. In particular,
$\text {pr}_U\inv (W)$ lies in a unique \ircompns~$C$ of $U\times _XX'$; and,
since $\text {pr}_U$ is a proper mapping, 
we must have $\text {pr}_U(C)=U$. Passing to a desingularization of~$C$, we 
get a 
commutative diagram of \holo mappings 
\begin{center}\begin{picture}(250,90)
\put(70,10){$U$}
\put(95,14){\vector(1,0){50}}
\put(155,10){$X$}
\put(115,0){$\Phi $}
\put(75,55){\vector(0,-1){30}}
\put(70,60){$U'$}
\put(60,35){$\alpha $}
\put(95,64){\vector(1,0){50}}
\put(155,60){$X'$}
\put(115,72){$\Phi '$}
\put(165,35){$\beta $}
\put(160,55){\vector(0,-1){30}}
\end{picture}
\end{center}
%
%
where $U'$ is a \con \cpx manifold of dimension~$m$,
$\alpha : U' @>>> U$ is a proper modification,
$S'\equiv \alpha\inv (S)=(\Phi ')\inv (R')$, and, if 
$W'=U'\setminus S'=\alpha\inv (W)$, then $\alpha $ maps the set 
$W'\setminus \alpha\inv (\sing U)$ biholomorphically 
onto~$\reg W$. We also get a \con 
covering space 
$\pi '=(\beta |_{M'})\inv \circ \pi : \Mtil @>>> M'$ and a 
lifting $\Phitil ' =
\Phitil \circ (\alpha | _{W'}) : W' @>>> \Mtil $ of
$\Phi '| _{W'}$. 

Therefore, if $L'=\beta ^*L$, then $\theta $ is a \holo $n$-form with values in
$\pi^*L=(\pi ')^*L'$ 
which is in~$L^2$ \wrt the metrics $\pi ^*h=(\pi ')^*\beta ^*h$ 
in $(\pi ')^*L'$ and 
$\pi ^*g=(\pi ')^*\beta ^*g$ on~$\Mtil $. Suppose the pullback of $\theta $ 
to~$W'$ extends to 
a section 
$$
\eta \in H^0(U', \ol ((\Phi')^*L')\otimes \Omega ^n_{U'}).
$$
Since $(\Phi')^*L'=\alpha ^*\Phi ^*L$ and $\alpha :U' @>>> U$ is a proper 
modification,
we have (by the definition of $\Omega ^n_U$ and Lemma~1.2)
$$
\alpha _*\bigl( \ol ((\Phi')^*L')\otimes \Omega ^n_{U'}\bigr)
=\ol (\Phi^*L)\otimes \Omega ^n_{U}.
$$
Hence $\eta $ determines an extension 
of $(\Phitil | _{\reg W})^*\theta $ to a section in
$H^0(U, \ol (\Phi^*L)\otimes \Omega ^n_{U})$.
Thus we may assume that $U$~and~$X$ are smooth and that $R$ is a divisor 
with normal
crossings in~$X$. In particular, $S=\Phi\inv (R)$ is a divisor in~$U$.

Since the lemma is entirely local, it suffices to extend the section near each 
point $x_0\in S$ and we may assume that $U=\Delta ^m$ is the unit polydisk centered 
at~$x_0=0$ in~$\C ^m$, that
$X=\Delta ^n$ is the unit polydisk centered at~$\Phi (x_0)=0$ in~$\C ^n$,
that $L$ is the trivial line bundle with the trivial metric on~$X$ 
(since all metrics are comparable
on \rel \cpt subsets), and that  $g$ is the restriction of the Euclidean metric
$g_{\C ^n}$ to~$M$ (since the $L^2$ condition on forms of type $(n,0)$ is 
independent of the choice of the metric on an $n$-dimensional manifold). 
We denote the 
coordinates in~$\C ^m$ by $z=(z_1,\dots , z_m)$, 
the 
coordinates in~$\C ^n$ by $\zeta =(\zeta _1,\dots , \zeta _n)$, and
the coordinate \fns of the mapping by $\Phi =(\Phi _1, \dots , \Phi _n)$.
Thus $\theta =fd(\zeta _1\circ \pi ) \wedge \dots \wedge 
d(\zeta _n\circ \pi )$ for some \holo
\fnns~$f$ which is square integrable on~$\Mtil $ with respect to 
$\pi ^*g_{\C ^n}$ and 
$\Phitil ^*\theta =(f\circ \Phitil)d\Phi _1\wedge \dots \wedge d\Phi _n$ on $W$. 
%
%
Since \holo sections extend past \anal sets of codimension at least~$2$,
we may assume that~$x_0\in \reg S$ and hence that~$S$ is 
the zero set of~$z_1$.  Since $R$ is a 
divisor with 
normal crossings, we may also assume that $R$ is the zero set of
$\zeta _1\cdots \zeta _k$. Finally, 
if $\setof {x\in S}{\Phi _j(x)=0}$ is nowhere dense in~$S$ for some~$j$, 
then, again, it suffices to 
consider a point~$x_0$ which avoids this zero set. Thus we may assume that
$$
S=\Phi \inv (R)=\{ \Phi _j=0 \} \text { for } j=1, \cdots , k.
$$

We now show that 
$(\Phi _1\cdots \Phi _k) \cdot (f\circ \Phitil )$ extends to a \holo 
\fn which vanishes along $S$. 
If~$x=(x_1,\dots ,x_m)$ is a point 
in~$W=U\setminus S=\Delta ^*\times \Delta ^{m-1}$ 
near~$x_0=0$ and $y=\Phi (x)=(y_1, \dots ,y_n)$, then we have
$r_j=|y_j|< 1/2$ for $j=1,\dots , n$ and $r_j>0$
for $j=1, \dots ,k$. Thus the polydisk
$$
P =\Delta (y_1;r_1)
\times \dots \times \Delta (y_k;r_k)
\times \Delta (y_{k+1};1/2)
\times  \Delta (y_n;1/2)
$$
centered at~$y$ is contained in $M=(\Delta ^*)^k\times \Delta ^{n-k}$
and is therefore evenly covered by $\pi : \Mtil @>>> M$. Hence~$\pi $
maps the \concompns~$\Ptil $ of~$\pi\inv (P)$ containing
$\ytil =\Phitil (x)$ isomorphically onto~$P$. The $L^{\infty }/L^2$-estimate 
now gives 
$$
|f(\ytil )|^2\leq (\Vol (\Ptil ))\inv 
\int _{\Ptil }
|f|^2 \, dV_{\pi ^* g_{\Cn }}.
$$
As $x$ approaches a point~$x_1$ in~$S$ near~$x_0$, $\Vol (\Ptil )=\Vol (P)$
will approach~$0$. Therefore, after multiplying both sides of the above 
inequality by $(r_1\cdots r_k )^2$ we get, since $|f|^2$ is integrable 
on~$\Mtil $, 
$$
|\Phi _1(x)\cdots \Phi _k(x) f(\Phitil (x))|^2=(r_1\cdots r_k)^2|f(\ytil )|^2
\leq 
\pi ^{-n}4^{(n-k)}\int _{\Ptil }
|f|^2 \, dV_{\pi ^* g_{\Cn }} @>>> 0 
$$
and the claim follows. 

For each $j=1,\dots ,k$, we have $\Phi _j=z_1^{\mu _j}h_j$ where 
$\mu _j =\text {ord}_S\Phi _j$  and $h_j$ is a unit. Therefore,
setting $\mu =\mu _1+\dots + \mu _k $ and 
$\psi =d\Phi _{k+1}\wedge \dots \wedge d\Phi _n$, we get 
\begin{align*}
d\Phi _1\wedge \dots \wedge d\Phi _n
&=z_1^\mu  dh_1\wedge  \dots \wedge dh_k \wedge \psi \\  
&\quad+z_1^{\mu -1} 
\sum _{j=1}^k \mu _j h_j (-1)^{j-1}dz_1 \wedge 
dh_1\wedge \dots \wedge \widehat {dh_j}\wedge 
\dots \wedge dh_k \wedge \psi .
\end{align*}
Since $z_1^{\mu }(f\circ \Phitil )$ extends to a \holo function 
which vanishes along~$S$, it follows that 
the $n$-form $\Phitil ^*\theta =(f\circ \Phitil )d\Phi _1\wedge \dots
\wedge d\Phi _n$
also extends holomorphically as claimed.  
\end{pf*}

\section{Results from the $L^2$~$\dbar $-method}

As described in Sect.~1, the proofs of the weak Lefschetz theorems rely
on a consequence of the $L^2$~$\dbar $-method 
(Andreotti-Vesentini~[AV],
H\"ormander~[Ho], Skoda~[Sk], Demailly~[D1]) 
which will be described in this section. 

Given a real-valued \fnns~$\vphi $ of class~$C^2$ on a \cpx manifold~$M$
of dimension~$n$, the {\it Levi form} $\lev \vphi $ of~$\vphi $ is the
Hermitian tensor defined by
$$
\lev \vphi = \sum _{i,j=1}^n
\frac {\partial ^2\vphi }{\partial z_i\partial \bar z_j} dz_id\bar z_j
$$
in local \holo coordinates $(z_1,\dots ,z_n)$. The function $\vphi $ is
said to be \plsh if $\lev \vphi \geq 0$ and \str \plsh if 
$\lev \vphi > 0$. If $(L,h)$ is a 
Hermitian \holo line bundle on a \cpx manifold $M$, then the 
{\it curvature tensor}~$\cal C (L,h)$ of~$(L,h)$ is given by
$$
\cal C(L,h)= \lev {-\log |s|^2}
$$
for any nonvanishing local \holo section~$s$ of~$L$. 
We will need the following
special case of a theorem of Demailly~[D1, Theorem~5.1]
concerning the $\dbar $-method for singular metrics with semi-positive
curvature. 
\begin{thm}[Demailly]  \it
Let $(E,h)$ be a Hermitian \holo line
bundle with semi-positive curvature (i.e.~$\cal C(L,h)\geq 0$) on a complete
K\"ahler manifold~$(M,g)$ of dimension~$n$. Suppose 
$\varphi : M @>>> [-\infty ,0]$
is a function which is of class~$\cinf $ outside a discrete 
subset~$S$ of~$M$ and, near each point $p\in S$, 
$\vphi (z)=A_p\log |z|^2$
where $A_p$ is a positive constant and $z=(z_1,\dots ,z_n)$ are local \holo
coordinates centered at~$p$. Assume that
$\cal C (E,he^{-\varphi })=\cal C (E,h)+\lev {\varphi }\geq 0$ 
on~$M\setminus S$
(and hence on~$M$ as the curvature of a singular metric) and let
$\lambda : M @>>> [0,1]$ be a \cont \fn \st 
$\cal C (E,h)+\lev {\varphi }\geq \lambda g$
on~$M\setminus S$. Then, for every $\cinf $~form~$\theta $ of type~$(n,1)$
with values in~$E$ on~$M$ which satisfies
$$
\dbar \theta =0 \quad \text {and} 
\quad \int _M\lambda \inv |\theta |^2_{h\otimes g^*}e^{-\varphi }
\, dV_g <\infty,
$$
there exists a $\cinf $ form~$\eta $ of type~$(n,0)$ with values in~$L$ on~$M$
\st
$$
\dbar \eta =\theta \quad \text {and} 
\quad \int _M |\eta |^2_{h\otimes g^*}e^{-\varphi }\, dV_g 
\leq \int _M\lambda \inv |\theta |^2_{h\otimes g^*}e^{-\varphi }
\, dV_g .
$$
\end{thm}

\begin{rem}
Demailly's theorem is much stronger than the above
special case. This special case also follows 
from Theorem~4.1 of~[D1], since one can approximate~$\varphi (z)$ by 
functions which locally have the form~$A_p\log (|z|^2+\epsilon )$ near the
nonsmooth points; or one can complete the metric on~$M\setminus S$.
\end{rem}

%
%

A well-known technique for producing sections with prescribed values on 
a discrete set gives the following:
\begin{thm}
Suppose $(L,h)$ is a Hermitian 
\holo line bundle on
an \irred reduced \cpx space~$X$ of dimension~$n$ and the curvature of~$h$ is
semipositive on~$X$ and positive at some point in~$X$. Then there exist a
positive integer~$\nu _0$ and a positive constant~$c_0$ which depend only 
on~$X$ and~$\cal C(L,h)$ and which have the following property. If~$\nu $ is an
integer with $\nu \geq \nu _0$, $R$ is a nowhere dense  \anal subset of~$X$ whose
complement~$M=X\setminus R$ is smooth and admits a complete K\"ahler metric, 
$(F,k)$ is a Hermitian \holo line bundle on~$X$ with semi-positive curvature,
$E_{\nu }=L^{\nu }\otimes F$, 
and $\pi : \Mtil @>>> M$ is a \con covering space of degree~$d$
($1\leq d\leq \infty $), then
$$
c_0\nu ^nd \leq 
\dim H^0_{L^2}(\Mtil , \ol (\pi ^* (E_{\nu } \otimes K_M))). 
$$
The $L^2$ condition is taken with respect to the Hermitian metric
$\pi ^*(h^{\nu }\otimes k)$ in $\pi ^*E_{\nu }$ 
and, for any choice of a Hermitian metric~$g$ on~$\Mtil $,
with respect to the Hermitian metric 
$g^*$~in~$K_{\Mtil }=\pi ^*K_M$ and~$g$ on~$\Mtil $
(the $L^2$-norm of an $(n,0)$-form does not depend on the choice of the
metric on the manifold).
\end{thm}

\begin{rems}
1. The curvature condition on~$L$
means that if $s$ is a nonvanishing  
\holo section of~$L$ on an open set~$W$, then $-\log |s|^2_h$ is \plsh
and, for some choice of~$W\neq \emptyset $, $-\log |s|^2_h$ is \str \plshns .

\noindent 2. The proof will also show that 
\begin{align*}
c_0\nu ^nd
&\leq c_0\nu ^n(d-1)+\dim H^0(X,\ol (E_{\nu })\otimes \Omega ^n_X) \\
&\leq 
\dim \biggl( 
H^0_{L^2}(\Mtil , \ol (\pi ^* (E_{\nu } \otimes K_M)))
+\pi ^*H^0(X,\ol (E_{\nu }) \otimes \Omega ^n_X)
\biggl);
\end{align*}
where the sum in the last expression takes place in 
$H^0(\Mtil ,\ol (\pi ^*(E_{\nu }\otimes K_M)))$.

\noindent 3. By a theorem of Demailly~[D1], $M=X\setminus R$ admits
a complete K\"ahler metric if, for example, $X$ is a complete
K\"ahler manifold and~$R$ is a \cpt \anal subset. In particular, any smooth
quasiprojective variety admits a complete K\"ahler metric. Thus we get as
a special case the following:
%
%
%
\end{rems}
\begin{cor}
Suppose $(L,h)$ is a positive Hermitian 
\holo line bundle on an \irred reduced projective variety~$X$
of dimension~$n$. Then there exist a
positive integer~$\nu _0$ which depends
on~$X$ and~$\cal C(L,h)$ and which has the following property. If~$\nu $ is an
integer with $\nu \geq \nu _0$, $R$ is a nowhere dense \anal subset
of~$X$ with smooth complement~$M=X\setminus R$, and 
$\pi : \Mtil @>>> M$ is a \con covering space of degree~$d$, then
$$
d\leq \dim H^0_{L^2}(\Mtil , \ol (\pi ^* (L^{\nu } \otimes K_M))).
$$
The $L^2$ condition is taken with respect to the Hermitian metric
$\pi ^*h^{\nu }$ in $\pi ^*L^{\nu }$ 
and, for any choice of a Hermitian metric~$g$ on~$\Mtil $,
with respect to the Hermitian metric 
$g^*$~in~$K_{\Mtil }=\pi ^*K_M$ and~$g$ on~$\Mtil $.
\end{cor}

\begin{pf*}{Proof of Theorem~2.2} By hypothesis, $\cal C(L,h) \geq 0$ on~$X$ and 
$\cal C(L,h)>0$ on
some \rel \cpt open subset~$W$ of~$X$. We may assume that $W\subset \reg X$
and that there exist
\holo coordinates $z=(z_1,\dots ,z_n)$ with $|z|<1/2$ on~$W$. 
Fix a nonempty \rel \cpt
open subset~$V$ of~$W$ and a~$\cinf $ \fnns~$\rho $ with \cpt support
in~$W$ \st $\rho \equiv 1$ on a \nbd of~$\overline V$, and, for
each point~$p\in V$, let~$\varphi _p$ be the~$\cinf $ \fn on~$X$ defined
by
$$
\varphi _p(x) =
\left\{ 
\begin{alignedat}{2}
&\rho (x)\log (|z(x)-z(p)|^2)& \quad  \text { if } x\in W \\
&0&\quad \text { if } x\in X\setminus W
\end{alignedat}
\right.
$$
Then $\supp \varphi _p =\supp \rho \subset W$,
$\varphi _p =\log (|z-z(p)|^2)$ (a \plsh \fnns ) on~$V$, and there is a 
positive constant~$a_0$ which does not depend on~$p$ \st
$a_0\cal C(L,h)+\lev {\varphi _p}$ is semipositive on~$X\setminus \{ p \} $ 
and 
positive on~$W\setminus \{ p \} $. Fix an integer~$\nu _0 > na_0$ 
(later, we will also choose $c_0$ to depend only
on~$a_0$ and~$n$). 

Let $\sing X\subset R\subset X$, $\pi : \Mtil @>>> M=X\setminus R$, and d be as in the
statement of the theorem and fix a point~$p$ in the nonempty open
set~$V\setminus R$. Given a multi-index 
$\alpha =(\alpha _1,\dots ,\alpha _n)\in \Z _{\geq 0}^n$,
a nonnegative integer~$\nu $, and a $\cinf $ section~$s$ of
$E_{\nu }\otimes K_X=L^{\nu }\otimes F\otimes K_X$ on a \nbd of~$p$, we
denote by $\partial ^{|\alpha |}s/\partial z^{\alpha }$ the corresponding
multiple derivative of~$s$ with respect to some fixed trivialization 
in~$L$~and~$F$ on a \nbd of~$p$ and the trivialization in~$K_X$ induced by the 
\holo coordinates $z=(z_1,\dots ,z_n)$ on~$W$. Similarly, if~$s$ is a~$\cinf $
section of $\pi ^* (E_{\nu }\otimes K_M)=\pi ^*E_{\nu } \otimes K_{\Mtil }$
on a \nbd of a point $q\in \pi\inv (p)$, then we denote by 
$\partial ^{|\alpha |}s/\partial z^{\alpha }$ the corresponding
multiple derivative of~$s$ with respect to the trivialization and local 
coordinates lifted from~$X$.

We will now apply the $\dbar $-method to show that if $\nu \geq \nu _0$, 
$\alpha =(\alpha _1,\dots ,\alpha _n)$ is a multi-index with 
$|\alpha |=\sum \alpha _j \leq (\nu /a_0)-n$, and $q\in \pi\inv (p)$, then
there exists a section
$$
s\in H^0_{L^2}(\Mtil , \ol (\pi ^* (E_{\nu } \otimes K_M)))
=H^0_{L^2}(\Mtil , \ol (\pi ^* E_{\nu }\otimes K_{\Mtil }))
$$
\stns , for every multi-index~$\beta $ with 
$|\beta | \leq (\nu /a_0)-n$ and for every point $r\in \pi\inv (p)$, we have
$$
\frac {\partial ^{|\beta |}s}{\partial z^{\alpha }}(r)=
\left\{ 
\begin{alignedat}{2}
&1& \quad  \text { if } \beta =\alpha \text { and } r=q \\
&0&\quad \text { otherwise}
\end{alignedat}
\right.
$$
By hypothesis, there exists a complete K\"ahler metric~$g$ on $M=X\setminus R$.
Let $\Etil _{\nu }=\pi ^*E_{\nu }$ for each~$\nu $, let $\htil =\pi ^*h$,
let $\ktil =\pi ^*k$, let $\gtil =\pi ^*g$, and let 
$\varphitil _p=\pi ^*\varphi _p$. We may choose a \rel \cpt \nbdns~$U$ 
of~$q$ in $\pi\inv (V) \setminus (\pi\inv (p)\setminus \{ q \} )$ and 
a~$\cinf $ section~$u$ of $\Etil _\nu \otimes K_{\Mtil }$ with \cpt
support in~$U$ \st $u$ is \holo on a \nbd of~$q$ in~$\Mtil $ and, for
every multi-index~$\beta $,
$$
\frac {\partial ^{|\beta |}u}{\partial z^{\alpha }}(q)=
\left\{ 
\begin{alignedat}{2}
&1& \quad  \text { if } \beta =\alpha \\
&0&\quad \text { if } \beta \neq \alpha 
\end{alignedat}
\right.
$$
Hence the form $\theta =\dbar u$ is a~$\cinf $~$\dbar $-closed 
$(n,1)$-form with values in~$\Etil _\nu $ and the support of~$\theta $
is a \cpt subset of $U\setminus \pi\inv (p)$ (since $\dbar u=0$ near~$q$).
By construction, there is also a \cont \fn $\lambda : \Mtil @>>> [0,1]$
\st $\lambda >0$ on~$\pi\inv (V)$ and
$$
\cal C(\Etil _\nu ,{\exp ({-\frac {\nu }{a_0}\varphitil _p})
\htil ^\nu \otimes \ktil })
= 
\nu \cal C(\pi ^*L,\htil ) + \cal C(\pi ^*F,\ktil )
+\frac {\nu }{a_0}\lev {\varphitil _p}
\geq \lambda \gtil 
$$
on~$\Mtil \setminus \pi\inv (p) $. 
Moreover, 
$$
\int _{\Mtil } \lambda \inv 
|\theta |^2_{\htil ^\nu \otimes \ktil \otimes  \gtil ^*}
e^{-\frac {\nu }{a_0}\varphitil _p} \, dV_{\gtil } < \infty ,
$$
because $\theta $ has \cpt support in 
$U\setminus \pi\inv (p)$, $\lambda >0$ on~$U$, and $\varphitil _p$ is
smooth on $U\setminus \{ q \}= U\setminus \pi\inv (p)$.
Applying Demailly's theorem (Theorem~2.1), one gets a $\cinf $ form~$\eta $
of type~$(n,0)$ with values in $\Etil _\nu $ on~$\Mtil $ \st
$$
\dbar \eta =\theta \quad \text {and} \quad
\int _{\Mtil } 
|\eta |^2_{\htil ^\nu \otimes \ktil \otimes  \gtil ^*}
e^{-\frac {\nu }{a_0}\varphitil _p} \, dV_{\gtil } < \infty .
$$
In particular, the $(n,0)$-form $s=u-\eta $ is in 
$H^0_{L^2}(\Mtil , \ol (\Etil_{\nu }\otimes K_{\Mtil }))$
because $\dbar s=0$ and $\varphitil _p\leq 0$.
Since $u$ is \holo near each point $r\in \pi\inv (p)$, so
is~$\eta $. Moreover, in suitable local \holo coordinates 
$w=(w_1,\dots ,w_n)$ centered at~$r$ in a \nbd~$Q$, we have
$\varphitil _p(w)=\log |w|^2$ and hence
$$
\int _Q|\eta |^2|w|^{-2\nu/a_0} \, dV<\infty ;
$$
where the notation for the metrics has been suppressed.
Therefore $\eta $ vanishes at~$r$ to an order greater than $(\nu /a_0)-n$.
%
%
Thus, if $\beta $ is a multi-index with $|\beta |\leq (\nu /a_0)-n$,
then $\partial ^{|\beta |}\eta /z^\beta $ vanishes at each point
$r\in \pi\inv (p)$ and the claim follows. 

The claim implies that if, for each $c\geq 0$, 
$b_c=\binom {[c]+n}{n}$ denotes the number of 
multi-indices~$\alpha $ satisfying $|\alpha | \leq c $, then we have
$$
\dim H^0_{L^2}(\Mtil , \ol (\Etil _{\nu }\otimes K_{\Mtil }))
\geq b_{(\nu /a_0)-n}\cdot d
$$
for each integer $\nu \geq \nu _0$. It is easy to see that 
$b_{(\nu /a_0)-n}\geq c_0\nu ^n$ for some positive constant~$c_0$ depending 
only on $a_0$~and~$n$ and the theorem now follows. 
\end{pf*}
\begin{rem} To obtain the inequalities given in the remark~2, 
we fix a point $r\in \pi\inv (p)$ and, for each point
$q\in \pi\inv (p)\setminus \{ r \} $ and each multi-index~$\alpha $
with $|\alpha |\leq (\nu /a_0)-n$, we form a section in 
$H^0_{L^2}(\Mtil , \ol (\Etil _{\nu }\otimes K_{\Mtil }))$
as in the above proof. We then get a collection of 
$b_{(\nu /a_0)-n}\cdot (d-1)$
linearly independent sections and the span of this collection
meets $\pi ^*H^0(X,\ol (E_\nu ) \otimes \Omega ^n_X)$ only in the zero 
section.
Therefore
\begin{align*}
c_0\nu ^n(d-1)&+\dim H^0(X,\ol (E_{\nu }) \otimes \Omega ^n_X) \\
&\leq 
\dim \biggl( 
H^0_{L^2}(\Mtil , \ol (\Etil _{\nu } \otimes K_{\Mtil }))
+
\pi ^*H^0(X,\ol (E_{\nu }) \otimes \Omega ^n_X)\biggl) .
\end{align*}
Finally, observe that if we take $\Mtil =M$ and 
$s\in H^0_{L^2}(M , \ol (E _{\nu }\otimes K_M))$, then, since the 
$L^2$~condition in the canonical bundle is independent of the choice of
the metric on the base manifold (provided one also takes the
associated metric in the canonical bundle), the 
pullback~$s'$
of $s$ to a desingularization~$X'$ of~$X$ is locally
in~$L^2$ with respect to a metric on~$X'$. Therefore $s'$ extends to
a section in~$H^0(X', \Omega ^n_{X'})$ and hence $s$ extends to
a section in~$H^0(X, \Omega ^n_X)$. It follows that
$$
c_0\nu ^n \leq \dim H^0(X,\ol (E_{\nu }) \otimes \Omega ^n_X).
$$
Thus we get all of the desired inequalities.
\end{rem}

\section{Further generalizations of the weak Lefschetz theorem}

Theorem~2.2 and the arguments given in the proofs of Theorem~1.1 
and Theorem~1.4 
now give the following generalization:
\begin{thm}
Let $U$ be an \irred reduced \cpx space, let $X$ be 
\con normal \cpx space of dimension $n>1$, 
let $\Phi :U @>>> X$ be a \holo mapping, and let 
$(L,h)$ be a Hermitian \holo line bundle on~$X$. Assume that 
\begin{enumerate}
\item[(i)] $\Phi (U)$ has nonempty interior, and 
\item[(ii)] The curvature of~$(L,h)$ is semipositive everywhere on~$X$ and 
positive at some point in~$X$.
\end{enumerate}
Then there exist a positive integer~$\nu _0$ and a positive constant~$c_0$ 
which depend only on~$X$ and (the curvature of)~$(L,h)$ 
\stns , if~$R$ is a nowhere dense  analytic subset of~$X$
whose complement~$\reg{X}\setminus R$ in $\reg{X}$ admits a complete K\"ahler 
metric, $V$~is a (nonempty) domain in~$U$, 
and~$\nu $~is an integer with $\nu \geq \nu _0$,
then 
\begin{align*}
c_0\nu ^n d  &\leq 
\dim H^0(V, \ol (\Phi ^*L^\nu )\otimes \Omega ^n_U); 
\end{align*} 
where $d$ is the index of the 
image~$G$ of $\pi_1(V\setminus \Phi\inv (R)) @>>> \pi _1(X\setminus R)$.
In particular, if 
$H^0 (V, \ol (\Phi ^*(L^\nu )) \otimes \Omega ^n_U)$
is finite dimensional for some choice of a sufficiently large~$\nu $,
then $G$ is of finite index. 
\end{thm}

%
%

Next, we show that for $U$ and $X$ smooth, there
exists a bound on the index in terms of the dimension of the
space of sections of an {\it invertible} sheaf. We first prove 
an elementary fact which relates sections of the pullback of the
canonical bundle to \holo $n$-forms (see also Sakai~[S]). 

\begin{lem}
Suppose 
$\Phi =(\Phi _1,\dots ,\Phi _n) : \Delta ^m @>>> \C ^n$ is a \holo mapping
and $\Phi _*$ has rank~$n$ at each point in $\Delta ^*\times\Delta ^{m-1}$.
We denote 
the coordinates in~$\C ^m$ and the coordinates in~$\Cn$ by 
$z=(z_1,\dots ,z_m)$
and
$\zeta =(\zeta _1,\dots ,\zeta _n)$, respectively. Let~$l\geq 0$ be the order
of vanishing of the \holo $n$-form $d\Phi _1\wedge \dots \wedge d\Phi _n$
along $\{ 0 \} \times \Delta ^{m-1}$. Then the mapping
$\ol (\Phi ^*K_{\C ^n}) @>>> \Omega ^n_{ \Delta ^m }$
given by
$$
s =f\Phi ^* (d\zeta _1\wedge \dots \wedge d\zeta _n) 
\mapsto fz_1^{-l}d\Phi _1\wedge \dots \wedge d\Phi _n
$$
maps $\ol (\Phi ^*K_{\C ^n}) $ isomorphically onto the sheaf of
\holo $n$-forms~$\theta $ \st 
$\theta _z\in \C (d\Phi _1\wedge \dots \wedge d\Phi _n)_z$ for each
point $z\in \Delta ^* \times \Delta ^{m-1}$ at which~$\theta $
is defined. 
\end{lem}
\begin{rem}
Here $\Phi ^* (d\zeta _1\wedge \dots \wedge d\zeta _n)$
denotes the pullback of $d\zeta _1\wedge \dots \wedge d\zeta _n$
as a section of $\Phi ^*K_{\Cn }$ while 
$d\Phi _1\wedge \dots \wedge d\Phi _n$ is the pullback as 
a form of type~$(n,0)$.  
\end{rem}

\begin{pf} Clearly, $z_1^{-l}d\Phi _1\wedge \dots \wedge d\Phi _n$ is 
a \holo $n$-form on~$\Delta ^m$, so we get an injective mapping
as described above. Conversely, suppose~$\theta $ is a \holo $n$-form on 
an open set $V\subset \Delta ^m$ and there exists a \holo \fnns~$h$ on
$V\cap (\Delta ^* \times \Delta ^{m-1})$ with 
$\theta =hd\Phi _1\wedge \dots \wedge d\Phi _n$ on 
$V\cap (\Delta ^* \times \Delta ^{m-1})$. We have
$$
\theta =\Sigma ' \theta _Idz_I \text { on } V \quad
\text {and} \quad
d\Phi _1\wedge \dots \wedge d\Phi _n
=\Sigma ' \beta  _Idz_I \text { on } \Delta ^m;
$$
where $\sum '$ denotes the sum over increasing multi-indices. In particular,
$$
l=\min _I (\text {ord} _{ \{ 0 \} \times \Delta ^{m-1} }\beta _I ) 
=\text {ord} _{ \{ 0 \} \times \Delta ^{m-1} }\beta _{I_0}
$$
for some multi-index~$I_0$, and, for each nonzero coefficient~$\beta _I$, 
we have $h=\theta _I/ \beta _I$ on 
$V\cap (\Delta ^* \times \Delta ^{m-1})$. 
Therefore $h$ is a 
meromorphic \fn on~$V$ with pole set contained in 
$\{ 0 \} \times \Delta ^{m-1}$ and 
$z_1^lh= \theta _{I_0}/ (\beta _{I_0}/z_1^l)$. Since the intersection of the
zero set of $\beta _{I_0}/z_1^l$ and $\{ 0 \} \times \Delta ^{m-1}$ 
has codimension at least~$2$ in~$\Delta ^m$ and since the
the pole set of $z_1^lh$ lies in this intersection, the pole set must 
be empty. Therefore $z_1^lh$
is \holo on~$V$ and the \holo section 
$s=z_1^lh\Phi ^*(d\zeta _1\wedge \dots \wedge d\zeta _n)$
of $\Phi ^*K_{\C ^n}$ maps to $\theta $
(on $V\cap (\Delta ^* \times \Delta ^{m-1})$ and hence on~$V$).
\end{pf} 

\begin{thm}
Let $U$~and~$X$ be \con \cpx manifolds of dimensions~$m$ and
$n>1$, respectively, let $\Phi :U @>>> X$ be a 
\holo mapping, and let~$(L,h)$ be a Hermitian \holo line bundle on~$X$.
Assume that   
\begin{enumerate}
\item[(i)] $\Phi $ has rank~$n$ at some point 
(i.e.~$\Phi (U)$ has nonempty interior), and 
\item[(ii)] The curvature of $(L,h)$ is semipositive everywhere on~$X$ and 
positive at some point in~$X$. 
\end{enumerate} 
Then 
there exist a positive integer~$\nu _0$ and a positive constant~$c_0$ which
depend only on~$X$ and (the curvature of)~$(L,h)$ and there exists 
an effective divisor~$D_0$ in~$U$ which depends only on the mapping
$\Phi : U @>>> X$ \stns , if~$R$ is a nowhere dense  analytic subset of~$X$
whose complement~$X\setminus R$ admits a complete K\"ahler metric,
$(F,k)$ is a Hermitian \holo line bundle on~$X$
with semipositive curvature, $V$~is a (nonempty) domain in~$U$, 
$\nu $~is an integer with $\nu \geq \nu _0$, $E_\nu =L^\nu \otimes F$, 
and $d$ is the index of the 
image~$G$ of $\pi_1(V\setminus \Phi\inv (R)) @>>> \pi _1(X\setminus R)$, 
then we have the estimates
\begin{align*}
c_0\nu ^nd
&\leq c_0\nu ^n (d-1) +\dim H^0 (X, \ol (E_\nu \otimes K_X)) \\
&\leq 
\dim H^0 (V, \ol (\Phi ^*(E_\nu \otimes K_X)\otimes [D_0])) \tag 1
\\
&\leq
\dim H^0(V, \ol (\Phi ^*E_\nu )\otimes \Omega ^n_U). 
\end{align*}   
\end{thm} 
\begin{rems}
1. If, for some positive integer~$k$, 
$L^k\otimes K\inv _X$ is semipositive,
then we may take $F=L^k\otimes K\inv _X$ and we get estimates which
do not involve the canonical bundle~$K_X$. 

\noindent 2. The divisor~$D_0$ which will be constructed is probably not
the optimal choice. 
\end{rems}

\begin{pf*}{Proof of Theorem~3.3}
Guided by Lemma~3.2, we first
describe~$D_0$.
The set 
$$
B=\setof {x\in U}{\text {rank}\, (\Phi _*)_x < n }
$$
is a nowhere dense  \anal subset of~$U$. Let $\{ A_i \} $ be the collection of all of
the \ircomps of~$B$ of dimension~$m-1$ whose image $\Phi (A_i)$ lies in some
nowhere dense  \anal subset of~$X$, let $A=\cup _i A_i$, and, for each~$i$, 
let~$l_i$ be the minimal order of vanishing along~$A_i$ of the 
$(n\times n)$-minor determinants of~$\Phi _*$. In other words, if 
$\Phi =(\Phi _1,\dots , \Phi _n)$ and 
$d\Phi _1\wedge \dots \wedge d\Phi _n=\sum 'a_Jdz_J$ with respect to local 
coordinates $(z_1,\dots ,z_m)$ near $x_0\in A_i$ in~$U$ and 
$(\zeta _1,\dots , \zeta _n)$ near~$\Phi (x_0)$ in~$X$, then
$$
l_i=\min _J (\text {ord}_{A_i} a_J).
$$
We define
$$
D_0=\sum l_iA_i.
$$

Given a \holo line bundle~$E$ on~$X$, 
Lemma~3.2 implies that we have an injective linear mapping
$$
H^0(U, \ol ((\Phi ^* (E\otimes K_X))\otimes [D_0]))
@>>>
H^0(U, \ol (\Phi ^* E)\otimes \Omega ^n_U)
$$
given as follows. Let~$t$ be a global defining section for~$[D_0]$
on~$U$. To each section 
$s\in H^0(U, \ol ((\Phi ^* (E\otimes K_X))\otimes [D_0]))$, we may
associate a \holo $n$-form~$\theta $ with values in 
$(\Phi ^* E)\otimes [D_0]$ and $\theta /t$ is a \holo $n$-form 
on~$U\setminus A$ with values in $\Phi ^* E$. But the lemma implies that,
near points of $A\setminus \sing B$, $\theta /t$ extends holomorphically
past~$A$. Thus $\theta /t$ extends to a 
\holo $n$-form on~$U\setminus (A\cap \sing B)$ with values in~$\Phi ^* E$. 
Since co$\dim \sing B \geq 2$, $\theta /t$ extends holomorphically to the
entire manifold~$U$. Thus we get a mapping $s \mapsto \theta /t$
(similarly, this mapping surjects onto the space of \holo $n$-forms 
with values in~$\Phi ^* E$ whose restriction to $U\setminus A$ comes
from a section of $\Phi ^*(E\otimes K_X)$). In particular, 
the third of the inequalities in~(1) holds.

Let 
\begin{center}\begin{picture}(250,80)
\put(5,10){$W=V\setminus S\subset V $}
\put(90,14){\vector(1,0){50}}
\put(145,10){$X \supset X\setminus R=M$}
\put(110,3){$\Phi $}
\put(125,55){\vector(3,-1){90}}
\put(110,60){$\Mtil $}
\put(17,25){\vector(3,1){90}}
\put(57,45){$\Phitil $}
\put(170,45){$\pi $}
\end{picture}
\end{center}
%
%
be a commutative diagram as in the proof of Theorem~1.1.  
We will show that if
$s\in H^0_{L^2}(\Mtil , \ol ((\pi ^*E_\nu  )\otimes K_{\Mtil }))$ 
for some~$\nu $ (with respect to metrics lifted 
from the base), then $(\Phitil ^*s)\otimes t$ extends to a unique 
\holo section of $(\Phi ^* (E_\nu \otimes K_X))\otimes [D_0]$ on~$V$.
By Lemma~1.3, the pullback of $s$ {\it as an $E_\nu $-valued \holo $n$-form}
extends to a $\Phi ^*E_\nu $-valued \holo $n$-form on~$V$. 
In particular, $\Phitil ^*s$ 
extends holomorphically as a section of 
$\Phi ^*(E_\nu  \otimes K_M)$ near each point at which
$\Phi _*$ is of maximal rank (by Lemma~3.2 with $l=0$). Moreover, 
$\Phitil ^*s$ extends holomorphically 
past \anal sets of codimension at least~$2$. Therefore, it suffices to
show that $(\Phitil ^*s)\otimes t$ extends holomorphically near each point
$x_0\in \reg S\cap \reg B$ at which $S\cap B$ is of dimension~$m-1$.
An \ircomp of~$B$ containing such a point~$x_0$ must also be an
\ircomp of~$S=\Phi \inv (R)$ and must therefore be one of the \ircompsns~$A_i$
of the support~$A$ of~$D_0$. Since the pullback of $s$ as an $E_\nu $-valued 
\holo $n$-form extends to~$V$, Lemma~3.2 and the definition of $D_0$ and $t$
now imply the claim. 

Clearly, if $s$ is a \holo section of $E_\nu \otimes K_X$ on $X$,
then $(\Phi ^*s)\otimes t$ is a \holo section whose restriction to~$V$ 
is an extension of $(\Phitil ^*\pi ^* s)\otimes t$. Thus we get an 
injective linear mapping of the subspace 
$$
\cal S=H^0_{L^2}(\Mtil , \ol ((\pi ^*E_\nu )\otimes K_{\Mtil }))+
\pi ^* H^0(X, \ol (E_\nu \otimes K_X))
$$ 
of 
$H^0(\Mtil , \ol ((\pi ^*E_\nu )\otimes K_{\Mtil }))$
into $H^0(V,\ol ((\Phi ^*(E_\nu \otimes K_X))\otimes [D_0]))$. 
We have, therefore 
\begin{align*}
\dim \cal S &\leq 
\dim H^0 (V, \ol (\Phi ^*(E_\nu \otimes K_X)\otimes [D_0])) \\
&\leq
\dim H^0(V, \ol (\Phi ^*E_\nu )\otimes \Omega ^n_U).
\end{align*}
The second remark following Theorem~2.2 now gives the inequalities~(1) for 
$\nu $ sufficiently large 
and for some constant~$c_0$ (both depending only on $(L,h)$ and $X$). 
\end{pf*} 

\begin{rems}
1. The proofs of Lemma~1.3 and Theorem~3.3 show that one can form a 
divisor $D_R$ 
which depends on~$R$, but which satisfies $D_R\leq D_0$ and gives
a sharper estimate for the index. For example, it suffices to include only
those \ircomps $A_i$ which are contained in $S=\Phi\inv (R)$, so one may choose
$D_R$ to have support contained in~$S$. Moreover, the proof of Lemma~1.3 shows
that if $\Phi (A_i)$  contains a point~$p\in R$ at which $R$ is a divisor with
normal crossings and $u$ is a defining \fn for $R$ near~$p$, then one 
may take the coefficient of~$A_i$ to be $-1+ \text {ord}_{A_i}(u\circ \Phi )$.
The proof also shows that, by choosing~$p$ so that this coefficient 
is minimal,  we get $D_R\leq D_0$.  

\noindent 2. Similarly, for $U$ a normal \nbd of a \cpt \cpx space~$Y$
and $X$ a smooth projective variety, one can 
find a uniform bound on the index
in terms of the dimension of a space of sections of a line
bundle pulled back from $X$ as in Theorem~1.1. More precisely, we have
the following:
\end{rems} 

\begin{cor} Let $\Phi :U @>>> X$ be a \holo mapping of a 
\con normal \cpx space~$U$ into a \con smooth projective variety~$X$
of dimension~$n>1$, let~$Y$ be a \con \cpt \anal 
subspace (not necessarily reduced) of~$U$, and 
let $\Uhat $ be the formal completion of~$U$ with respect to~$Y$.
Assume that  
\begin{enumerate}
\item[(i)] $\Phi (U)$ has nonempty interior, and 
\item[(ii)] 
$\dim H^0(\Uhat , \widehat {\ol (\Phi ^*L)})
<\infty $ for every \holo line bundle~$L$ on~$X$.
\end{enumerate}
Then there is a positive constant~$b$ depending only on the mapping
$\Phi :U @>>> X$ and the subspace~$Y\subset U$ \stns , if $R\subset X$
is a nowhere dense analytic subset of~$X$ and~$V$ is a \con \nbd of~$Y$ in~$U$,
then the image~$G$ of 
$\pi _1(V\setminus \Phi\inv (R)) @>>> \pi _1(X\setminus R)$
is of index at most~$b$ in~$\pi _1(X\setminus R)$.  
\end{cor}
%
%
\begin{pf*}{Sketch of the proof} First suppose $U$ is smooth and let
$D_0=\sum l_iA_i$ be the associated divisor in~$U$ as in the proof
of Theorem~3.3. By construction, each of the sets $\Phi (A_i)$ is contained
in some nowhere dense \anal subset of~$X$. By replacing $U$ by a \rel \cpt
\nbd of~$Y$, 
%
%
we may assume that there is a nowhere dense \anal subset~$C$ in~$X$ 
which contains all of these sets and that the collection of coefficients
$\seq li$ is bounded. Hence we may choose a positive \holo 
line bundle~$L$ on $X$ and a \holo section $t$ of~$L$ \st the 
divisor $D_1$ of the section $\Phi ^* t$ satisfies $D_1\geq D_0$.  

Now let $R$, $V$, and $G$ be as in the statement of the corollary
and let 
\begin{center}\begin{picture}(250,80)
\put(5,10){$W=V\setminus S\subset V $}
\put(90,14){\vector(1,0){50}}
\put(145,10){$X \supset X\setminus R=M$}
\put(110,3){$\Phi $}
\put(125,55){\vector(3,-1){90}}
\put(110,60){$\Mtil $}
\put(17,25){\vector(3,1){90}}
\put(57,45){$\Phitil $}
\put(170,45){$\pi $}
\end{picture}
\end{center}
be a commutative diagram as in the proof of Theorem~1.1.  
By the proof of Theorem~3.3 and the above remarks, 
if $s\in H^0_{L^2}(\Mtil ,\ol (\pi ^* L\otimes K_{\Mtil }))$,
then $(\Phitil ^*s)\otimes (\Phi ^* t)$ extends to a 
\holo section of $\Phi ^*(L^2\otimes K_X)$ on~$V$. 

If $U$ is \con and normal (but not necessarily smooth),
then we may form a desingularization $\alpha : U' @>>> U$ of~$U$
and a commutative diagram 
\begin{center}\begin{picture}(250,80)
\put(70,10){$U$}
\put(85,14){\vector(1,0){65}}
\put(155,10){$X$}
\put(115,0){$\Phi $}
\put(75,55){\vector(0,-1){30}}
\put(70,60){$U'$}
\put(60,35){$\alpha $}
\put(85,60){\vector(2,-1){70}}
\put(115,50){$\Phi '$}
\end{picture}
\end{center}
We may associate to $\Phi ' : U' @>>> X$ (after shrinking~$U$)
a line bundle~$L$ and
a section~$t$ as above, and we get the extension property for pullbacks
of $L^2$ sections as described. On the other hand, 
$U$ is normal, so
$\alpha _*\ol ((\Phi ')^*(L^2\otimes K_X))=\ol (\Phi ^*(L^2\otimes K_X))$.
Therefore the extension property also holds in~$U$, and
the usual argument now applies. 
\end{pf*}

We close this section by observing that
Theorem~3.1 has immediate consequences for pseudoconcave spaces.
An open subset $\Omega $ of a complex space~$X$ is said to have 
{\it pseudoconcave boundary in the sense of Andreotti}~[A]
if each point~$x_0\in \partial \Omega $ admits a fundamental system of
\nbdsns~$W$ in~$X$ \st $x_0$ is an interior point of
$$
\widehat {(W\cap \Omega )}_X
=\setof {x\in X}{|f(x)|\leq \sup _{W\cap \Omega }|f|\quad \forall \, 
f\in \ol (X) }.
$$
For example, by Proposition~10 of~[A],
if each \ircomp of~$X$ has dimension at least~$k>1$ and,
for each point~$x_0\in \partial \Omega $, there is a
$\cinf $~$(k-1)$-convex \fnns~$\vphi $ on a \nbdns~$W$ of~$x_0$ in~$X$
\st 
$$
\Omega \cap W =\setof {x\in W}{\vphi (x) >0},
$$
then $\Omega $ has pseudoconcave boundary in the sense of Andreotti.
A \con \cpx space~$X$ is said to be 
{\it pseudoconcave in the sense of Andreotti}~[A] if there
exists a nonempty \rel \cpt open subset~$\Omega $ which has 
pseudoconcave boundary in the sense of Andreotti and which meets each
\ircomp of~$X$. By a finiteness theorem of 
Andreotti~[A, Theorem~1], if $\cal F$ is a 
torsion-free coherent \anal sheaf on a locally \irred \con \cpx space~$X$
and $X$ is pseudoconcave in the sense of Andreotti, then
$\dim H^0(X,\cal F) < \infty $
(the case in which $X$ admits a $\cinf $~$(k-1)$-convex \exh 
\fnns, where the dimension of~$X$ is at least $k>1$ at each point, is due 
to Andreotti and Grauert~[AG]). 
Theorem 3.1 and Andreotti's finiteness theorem together give
the following:
\begin{cor}
Let $U$ be an \irred reduced \cpx space, 
let $X$ be a \con normal projective variety of dimension~$n>1$, and let 
$\Phi :U @>>> X$ be a \holo mapping. Assume that
$\Phi (U)$ has nonempty interior and that $U$ is pseudoconcave in the sense of 
Andreotti. 
Then there is a positive
constant~$b$ depending only on the mapping $\Phi : U @>>> X$ \stns , 
if~$Z$ is a nonempty Zariski open subset of~$X$, 
then the image of 
$\pi _1(\Phi\inv (Z)) @>>> \pi _1 (Z)$
is of index at most~$b$ in~$\pi _1(Z)$.
\end{cor}

\begin{rems}
1. Clearly, Theorem~3.1 also gives a version of the above theorem in which $X$
is not necessarily projective. 

2. There are many results concerning when
a \cpt \anal subset~$Y$ of an $m$-dimensional \cpx space~$U$
admits a strongly $(m-1)$-concave \nbdns ,
and hence when one may apply Andreotti's [A] (or Andreotti and 
Grauert's~[AG]) finiteness theorem as above.
For example, Okonek~[O] proved that
$Y$ admits a fundamental system of such \nbds if $N_{Y/U}$ 
is Finsler-$q$-positive, where $q=\dim Y$.
\end{rems}

\section{Burns' theorem}

The goal of this section is the following theorem:
\begin{thm}
Let $(X,g)$ be a \con complete Hermitian manifold
and let $M\subset X$  be a domain with
nonempty smooth \cpt boundary~$\partial M$ in~$X$. 
Assume that 
\begin{enumerate}
\item[(i)] $M$ is strongly pseudoconvex at each point of~$\partial M$;
\item[(ii)] There exists a Hermitian metric~$a$ in~$K_M$ and a 
constant $c>0$ \st $\cal C(K_M,a)\geq cg$ on~$M$; and
\item[(iii)] $X$ has dimension~$n\geq 3$. 
\end{enumerate} 
Then
$\Vol _g(M)<\infty $.
\end{thm}
\begin{rems}
1. Since $M$ admits a complete
K\"ahler metric, $\partial M$ is necessarily \con 
(see, for example, Proposition~4.4 below). 

\noindent 2. If, for example, the Ricci curvature of~$g$ is bounded
above by $-c$ on~$M$, then the associated metric $a=g^*$ in~$K_M$
satisfies the condition~(ii) since 
$$
\cal C(K_M,g^*)=-\text {Ric} \, (g) \geq cg.
$$


\noindent 3. Clearly, it is not necessary to assume that $M$ is a
domain in some larger manifold~$X$. The conclusion also holds if
$M=X$ and $M$ admits a $\cinf $ 
\fn which, along some end, is \str \plsh and exhaustive; since one can then 
replace $M$ by a suitable sublevel set of the function. It will, however,
be more convenient to have Theorem~4.1 stated for a domain as above. 

\noindent 4. As in the proofs of the weak Lefschetz theorems,
the idea is to apply finite dimensionality of a space of \holo 
sections of a line bundle to obtain a result about
the manifold. 
\end{rems} 

Theorem~4.1 and an analysis of the thick-thin decomposition as 
in~[BGS] together give as a conclusion the following theorem:
\begin{thm}[Burns~[B{]}]
Let $\Gamma $ be a torsion-free discrete group of 
automorphisms
of the unit ball~$B$ in~$\C ^n$ with $n \geq 3$
and let $M =\Gamma \setminus B$. Assume that the 
limit set~$\Lambda $ is a proper
subset of~$\partial B$ and that the quotient  
$\Gamma \setminus ((\partial B)\setminus \Lambda )$ 
has a compact component~$A$.
Then $M$ has only finitely many ends; all of which, except for the 
(unique) end corresponding to~$A$,
are cusps. In fact, $M$ is 
diffeomorphic to a \cpt manifold with boundary. 
\end{thm}


\begin{rem}
By applying a theorem of Lempert~[L] 
as in the proof of Theorem~4.1 below and the argument given by 
Siu and Yau~[SY], one can close up the cusps projectively.
In other words, 
$M\cong \Omega \setminus D$, where~$\Omega $ is a
strongly pseudoconvex domain in a smooth projective variety and~$D$ is
a (compact) divisor contained in~$\Omega $. The 
boundary component~$A$ 
corresponds to~$\partial \Omega $. 
\end{rem}

The main tool in the proof of Theorem~4.1 is 
Nadel and Tsuji's~[NT] $L^2$ version of 
Demailly's~[D2]
asymptotic Riemann-Roch inequality. 

\begin{thm}[Nadel-Tsuji~[NT{]}]
Suppose $(X,g)$ is a \con complete K\"ahler manifold
of dimension~$n$ and $(L,h)$ is a Hermitian \holo line bundle on~$X$
\st 
$$
\cal C(L,h)\geq cg
$$
for some constant $c>0$. Then
$$
\liminf _{\nu @>>> \infty } \nu ^{-n}
\dim H^0_{L^2}(X, \ol (K_X\otimes L^\nu )) \geq
\frac {1}{n!}\int _X \bigl( c_1(L,h)\bigr) ^n.
$$
\end{thm}
\begin{rems}
1. The Chern form $c_1(L,h)$ is the real 
form of type~$(1,1)$ (associated to the Hermitian tensor $\cal C(L,h)$)
given by 
$$
c_1(L,h)=-\frac {\sqrt {-1}}{2\pi }\partial \dbar \log |s|^2_h
$$
for any local nonvanishing \holo section~$s$ of~$L$. 

\noindent 2. As Nadel and Tsuji observed (see~[NT, Lemma~2.5]), 
if, in particular,  $X$ is pseudoconcave in the sense of 
Andreotti~[A] (see Sect.~3), 
then it follows that~$X$ has finite volume.

\noindent 3. The theorem is only stated in~[NT] for~$L$
the canonical bundle, but the proof of the general case is the same. 
The first point is that,
for a smooth \rel \cpt domain~$\Omega $ in~$X$ and for $\lambda >0$, 
one has Demailly's~[D2]
generalization of Weyl's asymptotic formula for the 
number of eigenvalues~$N_\Omega (\lambda )$ less than or equal 
to~$\nu \lambda $ for the Dirichlet problem for the Laplacian 
in
$K_X\otimes L^\nu $:
$$
\liminf _{\nu @>>> \infty } \nu ^{-n}
N_\Omega (\lambda )\geq
\frac {1}{n!}\int _\Omega  \bigl( c_1(L,h)\bigr) ^n.
$$
The second point is that for a $\cinf $ compactly supported form~$\alpha $
of type~$(n,1)$ with values in~$L^\nu $, the Bochner-Kodaira formula implies
that 
$$
\| \dbar \alpha \| _{L^2}^2+\| \dbar ^*\alpha \| _{L^2}^2 
\geq c\nu \| \alpha \| _{L^2}^2.
$$
With these slight changes in mind, the proof given in~[NT]
goes through. 
\end{rems}

We will also apply the following Hartogs type extension property:

\begin{prop} 
Let $(X,g)$ be a \con complete Hermitian manifold
of dimension~$n>1$ and let $M\subset X$ be a domain with nonempty
smooth \cpt strongly pseudoconvex boundary. Assume that
the restriction $g| _M$ of~$g$ to~$M$ is K\"ahler. 
Suppose~$f$ is a \holo \fn on $U\cap M$ for some
\nbdns~$U$ of~$\partial M$ in~$X$. 
Then there exists a \holo \fnns~$h$ on~$M$ \st $h=f$ near~$\partial M$. In
particular, $\partial M$ is \conns . 
\end{prop}
\begin{pf} We may assume that $M=\setof {x\in X}{\vphi (x) <0}$ for some
$\cinf $ \fnns~$\vphi $ on~$X$ which is \str \plsh on a \nbd 
of~$X\setminus M$ in~$X$. 
Since $g| _M$ is K\"ahler and $g$~is complete on~$X$, a theorem
of Nakano~[N] and of Demailly~[D1] implies that~$M$
admits a complete K\"ahler metric~$g'$. Moreover, the existence of~$\vphi $
implies that $(M,g')$ admits a positive Green's \fnns~$G$ which vanishes
along~$\partial M$. We normalize~$G$ so that, for each point $x_0\in M$,
$$
\lap _{\text {distr.}}G(\cdot , x_0) =-(2n-2)\sigma _{2n-1}\delta _{x_0};
$$
where $n=\dim X$, $\sigma _{2n-1}=\Vol (S^{2n-1})$, and 
$\delta _{x_0}$ is the Dirac \fn at~$x_0$. 

Fix a $\cinf $ \fnns~$\lambda $ with \cpt support in~$U$ \st 
$\lambda \equiv 1$ on a \nbd of~$\partial M$ and let~$\alpha $ be the 
$\dbar $-closed compactly supported form of type~$(0,1)$ on~$M$ given by
$\alpha =\dbar (\lambda f)$ (extended by~$0$ to~$M$). Then the \fnns~$\beta $
defined by
$$
\beta (x) =-\frac {1}{(2n-2)\sigma _{2n-1}}
\int _M G(x,y) \dbar ^* \alpha (y)\, dV_{g'}(y)
$$
is a $\cinf $ bounded \fn with finite energy 
(i.e.~$\int _M |\nabla \beta | ^2 \, dV_{g'} <\infty $), 
$\lap \beta = \dbar ^* \alpha $, and $\beta $ vanishes on~$\partial M$. Hence
$\gamma \equiv\alpha -\dbar \beta $ 
is an $L^2$ harmonic form of type~$(0,1)$ and
the Gaffney theorem~[G] implies that~$\gamma $ is closed
(and coclosed). In particular, $\bar \gamma $ is a \holo $1$-form on~$M$ and
$\beta $~is \plh on $W\cap M$ for some \nbdns~$W$ 
of~$X\setminus M$ in~$X$.

We will show that~$\beta $ vanishes near~$\partial M$. Fix $a<0$ so
close to~$0$ that~$\vphi $ is \str \plsh on 
$V=\setof {x\in M}{\vphi (x) >a}$ and $V\subset\subset W$. If 
$\rho $ is the real part or the imaginary part of~$\beta $ and $\rho $
does not vanish identically near~$\partial M$, then we may choose
a nonzero regular value~$b$ of~$\rho $ contained in~$\rho (V)$.
Since $b\neq 0$ and $\rho $ vanishes on~$\partial M$, $\rho\inv (b)$
avoids~$\partial M$. Thus the restriction of~$\vphi $ to~$\rho\inv (b)$
assumes its maximum at some point~$x_0\in V\subset W\cap M$ (with $\vphi (x_0)>a$). 
But the leaf~$L$ through~$x_0$ of the foliation determined by the 
\holo $1$-form~$\partial \rho $ on~$V\cap M$ is contained in~$\rho\inv (b)$,
so $\vphi | _L$ also assumes its maximum at~$x_0$.
Since~$\vphi $ is \str \plsh on~$V$, we have arrived at a contradiction.
Therefore $\beta $ vanishes near~$\partial M$. Hence 
$\gamma =\alpha -\dbar \beta $ vanishes near~$\partial M$ and, therefore,
on all of~$M$, since~$\bar \gamma $ is a \holo $1$-form. Thus the
\fnns~$h\equiv \lambda f-\beta $ is \holo on~$M$ (since $\dbar h=\gamma =0$)
and equal to~$f$ near~$\partial M$. In particular, since one can
take $f$ to be a 
locally constant \fn which separates distinct components of~$\partial M$, 
$\partial M$ is \conns . 
\end{pf}

\begin{pf*}{Proof of Theorem~4.1} Since $n\geq 3$, one can apply a
theorem of Rossi~[R] to 
``fill in the holes'' and obtain a \con Stein space~$Y$ with isolated
singularities, a \rel \cpt pseudoconvex domain~$N$ in~$Y$ 
containing~$\sing Y$, and a bi\holo mapping $\Phi : U @>>> V$ of a \nbdns~$U$
of~$\partial M$ in~$X$ onto a \nbdns~$V$ of~$\partial N$ in~$Y$ \st
$\Phi (U\cap M)=V\cap N$. Since $N$ may be embedded into a Euclidean space,
Proposition~4.4 implies that~$\Phi $ extends to a \holo
mapping $M\cup U @>>> Y$, which we also denote by~$\Phi $,
and $\Phi (M)\subset N$. Next, by a theorem
of Lempert~[L], one can form a ``cap'' on~$N$.
That is, we may assume that~$Y$ is an affine algebraic variety.
By forming the closure~$\overline {Y}$ of~$Y$ in a projective space and
desingularizing~$\overline {Y}$ at infinity, we get a projective 
variety~$Z$ with isolated singularities \st 
$\sing Z\subset N\cup V \subset Z$. Finally, by replacing~$X$ by
$$
(M\cup U) \cup (V\cup (Z\setminus \overline N))\bigg/  
x\in U \sim \Phi (x) \in V
$$
and by replacing the metric $g$ by any extension of $g| _M$ to the
new manifold, we may assume that we have a \holo mapping
$\Phi : X @>>> Z$ \st $\Phi (M) \subset N$ and $\Phi $ maps 
$(X\setminus M)\cup U$ biholomorphically onto $(Z\setminus N)\cup V$.
In particular, since $X\setminus \overline M \subset\subset X$, it follows
that $X$ is pseudoconcave in the sense of Andreotti. 

Now let~$H$ be a positive Hermitian \holo line bundle on~$Z$. Then $\Phi ^*H$
is semipositive on~$X$ and positive on $(X\setminus M)\cup U$. On the
other hand, by shrinking~$M$ slightly and extending the Hermitian 
metric~$a$, we may assume that $K_X$ admits a Hermitian metric whose
curvature is greater than or equal to $cg$ at each point of~$M$. 
It follows that if~$m$ is a
sufficiently large positive integer and 
$L=K_X\otimes \Phi ^*H^m$,
then~$L$ admits a Hermitian metric~$h$ \st 
$\cal C(L,h)\geq cg$
on~$X$. In particular, $g'=\cal C(L,h)$ is a complete K\"ahler metric on~$X$.
Therefore, by the $L^2$ Riemann-Roch inequality of Nadel and 
Tsuji (Theorem~4.3), we have, for every sufficiently large positive
integer~$\nu $,
$$
1+ \, \nu ^{-n}
\dim H^0_{L^2}(X, \ol (K_X\otimes L^\nu )) \geq
\frac {1}{n!}\int _X \bigl( c_1(L,h)\bigr) ^n
\geq c^n\pi ^{-n} \int _X \, dV_{g};
$$
(where the Hermitian metric in $K_X\otimes L^\nu $ is
$(g')^*\otimes h^\nu $).  Since,  by
Andreotti's finiteness theorem~[A]
(or by~[AG]), the left-hand side is finite, we get
$$
\Vol _g(M)\leq \Vol _g(X) <\infty .
$$
\end{pf*}

\begin{rem}
By a version of the $L^2$~Riemann-Roch
inequality due to Takayama~[T], it is only necessary to assume
in the hypothesis~(ii) that $\cal C(K_M,a)\geq cg$ outside a \rel 
\cpt \nbd of $\partial M$ in $X$. 
\end{rem}

%
%

\bibliographystyle{amsplain}

\end{document}